\documentclass[12pt]{article}
\usepackage{latexsym}
\usepackage{amsfonts}
\usepackage{amssymb}
\usepackage{amsmath}
\newtheorem{theorem}{Theorem}
\newtheorem{corollary}{Corollary}
\newtheorem{proposition}{Proposition}
\newtheorem{remark}{Remark}
\DeclareMathOperator*{\tends}{\longrightarrow}

\begin{document}
\author{Gheorghe Minea}
\title{Entropy conditions for quasilinear first order equations on nonlinear fiber bundles\\
with special emphasis on the equation of 2D flat projective structure. I.}
\maketitle
\begin{titlepage}
\begin{abstract}
Taking only the characteristics as absolute, in the spirit of Arnold [1], we give an independent of coordinates formulation of general
variational entropy inequalities for quasilinear equations of first order, that locally read as Kruzhkov inequalities, in terms 
of certain ``entropy  densities'', and in the case of the equation of 2D flat projective structure we get the expression of the general 
entropy density from its abstract Rankine-Hugoniot rule for shocks using the projective geometry of the plane.
\end{abstract} 
\end{titlepage}
\section{Introduction}

Our concern here is a local and independent of coordinates formulation, in the space-time continuum, of entropy conditions of Kruzhkov 
type for quasilinear first order equations, aiming at a geometric understanding of the meaning of time and unique determination in the 
future of the entropy solutions with shocks. The need for this inspection comes from three observations: first, that the entropy 
condition is local as the equation itself; next, that it is not invariant at the local, nonlinear in fibers, bundle transformations 
that, however, map graphs of local classical solutions into graphs of local classical solutions; and third, that there are nonlinear 
in fibers transformations of that kind that induce changes of coordinates in the space-time continuum mixing space and time coordinates.
In order to keep track of the entropy condition we have then to consider the solutions as sections of a nonlinear fiber bundle over a 
base where a priori there is not a distinguished time coordinate.\\
For instance, we prove here that a local, nonlinear in fibers, bundle transformation maps graphs of local classical solutions of the 
equation of 2D flat projective structure into surfaces of the same kind if and only if it is the natural lift of a projective 
transformation of the (space-time) plane to the projectivized of the tangent bundle to it; thus the equation of 2D flat projective 
structure may be naturally understood as an equation for sections of the (nontrivial) projectivized tangent bundle to the projective 
plane (whose fibers are projective lines, i.e. cercles). Along these lines, using  E. Cartan's theory of projective curvature, we 
succeed to characterize the quasilinear equations that can be locally reduced by a bundle change of coordinates to the equation of
2D flat projective structure.\\
The manifold of 1-jets of sections of a fiber bundle of 1-dimensional (nonlinear) fiber bears a canonical contact structure so that the 
theory of characteristics from [1] applies also to nonlinear first order partial differential equations for such sections. Moreover,
the set of 1-jets of sections of given value at a given point has a canonical affine structure so that the quasilinear equations for the 
sections of these general fiber bundles have a well defined meaning.\\
We show here that in this framework the entropy condition is caught by a section, that we call ``entropy density'', of a certain line 
bundle, over the total space of the nonlinear bundle, derived from two other line bundles over the same base; and in fact the ``entropy 
density'' is, for each point of the total space of the nonlinear bundle, a non-zero odd 1-form (in the sense of de Rham) from the tangent 
to the fiber to the tangent to the (reduced) characteristic at the same point. It then defines both the characteristic directions and a way 
(or orientation) on each characteristic: the way of time flowing.\\
In the special case of the equation of 2D flat projective structure the Rankine-Hugoniot rule of shock formation may be seen as taking 
the barycenter of the jump interval in the nonlinear fiber with respect to the measure on that fiber defined by the entropy density. We 
characterize those  barycentric maps in terms of the projective geometry of the line and of its tautologic fiber bundle and prove that 
the entropy density is uniquely determined (up to a positive factor) by its barycentric map. We succeed in this way to write down the 
general variational entropy inequality for the equation of 2D flat projective structure only in terms of its abstract Rankine-Hugoniot 
rule for shocks.

\section{Entropy density and Kruzhkov type\\
variational inequality on\\
nonlinear fiber bundles}

\S 2.1 \texttt{Some intrinsic operations with odd differential forms}

The geometric formulation of the variational inequality involves operations with odd differential forms; we recall the definition from 
de Rham [5], using however other notations. For a real vector space $V$ of $\dim V=d$ we denote $V^{\wedge p}$ the $p$-th exterior power 
of $V$,
\begin{equation}\label{1}
S(V)=\{ s:V^{\wedge d}\longrightarrow \textbf{R} \: \arrowvert \:s(\lambda v)=\mathrm{sgn}(\lambda)\cdot s(v),\forall\lambda\in
\textbf{R},\forall v\in V^{\wedge d} \}
\end{equation} 
the 1-dimensional space of odd 0-forms, or odd scalars, on $V$,
\begin{equation}\label{2}
(V^{\wedge p})^{\ast}\otimes S(V),\;0 \leq p \leq d,
\end{equation} 
the space of odd p-forms on $V$ and
\begin{equation}\label{3}
 \Omega(V)=\{\rho:V^{\wedge d}\longrightarrow \textbf{R} \: \arrowvert \:\rho(\lambda v)=|\lambda|\cdot \rho(v),\forall\lambda\in
\textbf{R},\forall v\in V^{\wedge d} \}
\end{equation} 
the space of densities on $V$. We consider also the set of orientations of $V$ :\\
\begin{equation}\label{4}
O(V)=\{s\in S(V)\arrowvert\;|s(v)|=1,\;\forall v\in\;V^{\wedge d}\smallsetminus \{0\}\}.
\end{equation} 
For $\alpha\in(V^{\wedge d})^{\ast}$ it is well defined
\begin{equation}\label{5}
 |\alpha|\in \Omega(V),
\end{equation}
and if $\alpha\neq 0$
\begin{equation}\label{6}
 \mathrm{sgn}(\alpha)\in O(V).
\end{equation}
We have a canonical isomorphism
\begin{equation}\label{7}
\Phi : \Omega(V) \widetilde{\longrightarrow}(V^{\wedge d})^{\ast}\otimes S(V)
\end{equation} 
given by
\begin{equation}\label{8} 
(\Phi\rho)(v_{1}\wedge\dots\wedge v_{d})(w_{1}\wedge\dots\wedge w_{d})=\rho(v_{1}\wedge\dots\wedge v_{d})\cdot\mathrm{sgn}(v_{1}
\wedge\dots\wedge v_{d}:w_{1}\wedge\dots\wedge w_{d})
\end{equation} 
where $v:w\in\textbf{R}$, for $v, w$ in the same 1-dimensional space, with $w\neq 0$, has the meaning that $(v:w)\cdot w=v$. For $W
\subset V$ subspace of $\dim W=p$, with $\pi:V\rightarrow V/W$ and $\iota:W\rightarrow V$ canonical, the mapping
\begin{equation}\label{9} 
(f_{1}\wedge\dots\wedge f_{p})\otimes(\pi g_{1}\wedge\dots\wedge \pi g_{d-p})\longmapsto \iota f_{1}\wedge\dots\wedge\iota f_{p} \wedge 
g_{1}\wedge\dots\wedge g_{d-p} 
\end{equation} 
is well defined and establishes a canonical isomorphism
\begin{equation}\label{10}
W^{\wedge p}\otimes (V/W)^{\wedge (d-p)} \widetilde{\longrightarrow} V^{\wedge d},
\end{equation} 
wherefrom the interesting for us isomorphism
\begin{equation}\label{11}
\Omega(W)\otimes \Omega(V/W)\widetilde{\longrightarrow}\Omega(V).
\end{equation}
We denote the image of this mapping as $\rho \otimes \sigma \mapsto\rho \rtimes\sigma $ so that in accordance with the formula (\ref{9}) 
we have.:
\begin{equation}\label{12} 
(\rho \rtimes\sigma)(v_{1}\wedge\dots\wedge v_{d})= \rho(Pv_{1}\wedge\dots\wedge Pv_{p})\sigma(\pi v_{p+1} \wedge\dots\wedge\pi v_{d})        
\end{equation}
if $\pi v_{p+1} \wedge\dots\wedge\pi v_{d}\neq 0$ and $P: V\rightarrow W$ is the projection corresponding to the direct sum decomposition 
$$V=W\dotplus\sum_{j=p+1}^{d}\textbf{R}\cdot v_{j}.$$ Next for $\mu\in\Omega(V)$ and $X\in V$ the contraction $i_{X} \mu\in (V^{\wedge 
(d-1)})^{\ast}\otimes S(V)$ is defined as the odd $(d-1)$-form
\begin{equation}\label{13}
 i_{X} \mu (v_{1}\wedge\dots\wedge v_{d-1})(w_{1}\wedge\dots\wedge w_{d})=(\Phi\mu)(X\wedge v_{1}\wedge\dots\wedge v_{d-1})(w_{1}\wedge
\dots\wedge w_{d})
\end{equation} 
(see the definition (\ref{8})). Now for $\rho\in\Omega(W)\otimes V,\;\phi\in\Omega(V/W)$ we define the contracted tensor product
\begin{equation}\label{14}  
\rho\rightthreetimes\phi=i_{X}(\lambda\rtimes\phi),
\end{equation} 
if $\rho=\lambda\otimes X,\;\lambda\in\Omega(W),\;X\in V$; the existence of $\lambda$ and $X$ comes from $\dim \Omega(W)=1$ and they 
are unique up to a non-zero factor that does not affect the definition. Of course $\rho\rightthreetimes\phi\in (V^{\wedge (d-1)})^{\ast}
\otimes S(V)$.\\
If $M$ is a differential manifold and $H\rightarrow M$ is a vector bundle over $M$ we denote by $\textsl{C}^{k}\Gamma(H)$ the space of 
its $\textsl{C}^{k}$ sections and by $\textsl{C}_{0}^{k}\Gamma(H)$ the space of sections of compact support. We identify a density 
$\rho\in\textsl{C}\Gamma(\Omega(TM))$ with the measure it defines and write $$\int \phi(z)\rho(\mathrm{d}z)$$ for the integral of $\phi$ 
with respect to it.\\
We will call $nonlinear\;fiber\;bundle$ a surjective submersion $\pi:F\rightarrow M$; in that case the tangent to the fiber 
$F_{z}=\pi^{-1}(\{z\})$ in $f$ (with $\pi(f)=z$) becomes a locally trivial vector bundle $T^{0}F\rightarrow F$ with 
$T^{0}_{f}F=:T_{f}F_{\pi(f)}=\mathrm{ker}\;T_{f}\pi$, and with the quotient vector bundle $TF/T^{0}F$; if $\pi^{\ast}(TM)\rightarrow F$ 
denotes the usual inverse image of the vector bundle $TM\rightarrow M$ through the map $\pi: F\rightarrow M$, the natural isomorphism of 
vector bundles\\
$$TF/T^{0}F\widetilde{\longrightarrow} \pi^{\ast}(TM)$$\\
will be understood.\\
If $\lambda\in\textsl{C}\Gamma(\Omega(T^{0}F))$, then $\lambda^{z}=:\lambda\arrowvert _{F_{z}}$ defines a measure over $F_{z}$. On the other 
hand, for $\phi\in\textsl{C}\Gamma(\Omega(TF/T^{0}F))$ the restriction $\phi^{z}=:\phi\arrowvert_{F_{z}}$ defines a vector function on 
$F_{z}$ with values in $\Omega(T_{z}M)$. Then \\
$$\int_{F_{z}} \psi(f)\phi^{z}(f)\lambda^{z}(\mathrm{d}f)\in \Omega(T_{z}M)$$\\
and the following equality
\begin{equation}\label{15}
\int_{F} \psi(f)(\lambda\rtimes\phi)(\mathrm{d}f)=
\int_{M}\left(\vphantom{\int_{F_{z}}}\int_{F_{z}}\psi(f)\phi^{z}(f)\lambda^{z}(\mathrm{d}f)\right)(\mathrm{d}z)
\end{equation} 
gives the interpretation, in terms of measures, of the product defined in fibers according to (\ref{12}).\\
Similarly, for $\rho\in\textsl{C}\Gamma(\Omega(T^{0}F)\otimes TF)$ and $\phi\in\textsl{C}\Gamma(\Omega(TF/T^{0}F))$ it is defined, through 
the contracted product (\ref{14}) taken in fibers, the odd (d-1)-form on $F$, if $d$ denotes $\dim F$:
\begin{equation}\label{16}
\rho\rightthreetimes\phi\in\textsl{C}\Gamma((T^{\ast}F)^{\wedge(d-1)}\otimes S(TF)).
\end{equation}
 
\S2.2 \texttt{The generalized function $I(\rho,\sigma,G)$ and its local disintegration}

From now on we consider only the case when $\dim F_{x}=1$ for $x\in M$.\\
Let $\rho\in\textsl{C}^{1}\Gamma(\Omega(T^{0}F)\otimes TF)$ be such that
\begin{equation}\label{17}
T_{f}\pi\;\rho_{f}\neq 0,\;\;\forall\; f\in F, 
\end{equation}
and define its \textit{characteristic directions} by
\begin{equation}\label{18}
D_{f}=\rho_{f}(T^{0}_{f}F)-\rho_{f}(T^{0}_{f}F),
\end{equation}
i.e. the 1D sub-space spanned by the range of $\rho_{f}$ in $T_{f}F$ (which is a half-line in $D_{f}$). If we start from the smooth 
1D sub-bundle $D$ of $TF$, with $T_{f}\pi\;D_{f}\neq 0,\;\forall\; f\in F$, we think 
\begin{equation}\label{19}
\rho\in\textsl{C}^{1}\Gamma(\Omega(T^{0}F)\otimes D);
\end{equation}
in this case, apart from $D$, it is enough to know $T\pi\;\rho\in\textsl{C}^{1}\Gamma(\Omega(T^{0}F)\otimes T\pi\;D)$,
where $(T\pi\;D)_{f}=:T_{f}\pi\;D_{f}$ is a sub-bundle of $\pi^{\ast}(TM)$ and $(T\pi\;\rho)_{f}=T_{f}\pi\;\rho_{f}$.\\
Let $\sigma\in\Gamma(F\arrowvert U)$ be a Lebesgue measurable section of $\pi:F\rightarrow M$ for which 
$\exists G=\mathring{G}\subseteq F$ such that $\sigma(U)\subseteq G,\;\pi(G)=U$ with the properties
\begin{equation}\label{20}
 G_{z}:=\pi^{-1}(\{z\}\cap G\;noncompact\;connected\;\forall\;z\in\pi(G)
\end{equation}
and $\forall\;z_{0}\in U\; \exists V=\mathring{V}\subseteq U,\;V\ni z_{0},
\;\exists \;\sigma_{1},\;\sigma_{2}\in\textsl{C}\Gamma(G\arrowvert V)$ such that
\begin{equation}\label{21}
\sigma_{z}\in \lvert\sigma_{1z},\sigma_{2z}\lvert_{G_{z}},\;z\in V,
\end{equation}
where $\lvert a,b\lvert_{G_{z}}$ denotes the closed nonoriented interval included in the open arc $G_{z}$ of ends $a,\;b\in G_{z}$. 
Here $\sigma_{1},\;\sigma_{2}$ are continuous, while $\sigma$ may be not; however, the set of limit points of $\sigma$ in $z$ is 
contained in the interval from the right hand side of (\ref{21}). We may call \textit{G open layer bounding} $\sigma$. The meaning of these 
definitions is made clear by the following\\  

\begin{proposition} $In\; the\; hypothesis\; (\ref{20})\; on\; G, \;\forall z_{0}\in U\;\exists V=\mathring{V}\subseteq U,\; V\ni z_{0},\\
\exists \Phi:\pi^{-1}(V)\cap G \rightarrow \textbf{R}^{m}\times\textbf{R}\;diffeomorphism\; on\; an\; open\; subset,\; mapping\; fibers\;\\
into\; fibers\; (i.e.\; a\; bundle\; map).$ \\
\end{proposition}

Next we consider the covering
\begin{equation}\label{22}
\varpi:\tilde{G}\rightarrow G
\end{equation}
where the fiber is defined by 
\begin{equation}\label{23}
 \tilde{G}_{g}=\varpi^{-1}(\{g\})=O(T_{g}G_{\pi(g)})
\end{equation}
i.e. the set of orientations of the space $T_{g}G_{\pi(g)}=\mathrm{ker}\;T_{g}\pi$ (see (\ref{4})). We will consider also the 
nonlinear fiber bundle\\
\begin{equation}\label{24}
 \tilde{\pi}:\tilde{G}\rightarrow U,\;\;\tilde{\pi}:=\pi\circ\varpi.
\end{equation}
The space of test functions for $I(\varrho,\sigma,G)$ will be $\textsl{C}_{0}^{\infty}\Gamma(\Omega(TG))$; the elements of that 
functional space being smooth densities, we chose the term of generalized function on $G$ for an element of its dual.\\
For a density $\psi\in\textsl{C}_{0}^{\infty}\Gamma(\Omega(TG))$ we define its $fiber\;primitive$
\begin{equation}\label{25}
 \int\psi\in\textsl{C}^{\infty}\Gamma(\Omega(T\tilde{G}/T^{0}\tilde{G}))
\end{equation}
as follows: for $\varepsilon\in\tilde{G},\;\varepsilon=(g,\omega),\;g=\varpi(\varepsilon),\;\omega\in O(T_{g}G_{\pi(g)})$ defines an 
orientation on the arc $G_{\pi(g)}$ so that it is well defined the interval
\begin{equation}\label{26}
 (-\infty,\;g)_{\omega}=\{h\in G_{\pi(g)}\arrowvert \;h<g \;with\; respect\; to\; \omega\};
\end{equation}
next, for $g\in G_{z}$, according to (\ref{11}), $\Omega(T_{g}G)\widetilde{\longrightarrow}\Omega(T_{g}G_{z})\otimes\Omega(T_{z}U)$, 
so that\\
$\psi_{g}\in\Omega(T_{g}G_{z})\otimes\Omega(T_{z}U),\;\forall g\in G_{z}$. Thus it is well defined 
\begin{equation}\label{27}
(\int\psi)_{\varepsilon}= \int_{(-\infty,\;g)_{\omega}} \psi\;
\in\Omega(T_{z}U)\widetilde{\longrightarrow}\Omega(T_{\varepsilon}\tilde{G}/T_{\varepsilon}^{0}\tilde{G}).
\end{equation}
In view of the Proposition 1, $\int\psi \in \textsl{C}^{\infty}\Gamma(\Omega(T\tilde{G}/T^{0}\tilde{G}))$ for
$\psi\in\textsl{C}_{0}^{\infty}\Gamma(\Omega(TG))$.\\
As $\varpi:\tilde{G}\rightarrow G$ is a local diffeomorphism and a bundle map between fiber bundles over $U$, for  
$\rho\in\textsl{C}^{1}\Gamma(\Omega(T^{0}G)\otimes TG)$ it is well defined 
\begin{equation}\label{28}
 \varpi^{\ast}\rho\in\textsl{C}^{1}\Gamma(\Omega(T^{0}\tilde{G})\otimes T\tilde{G})
\end{equation}
by the relation
\begin{equation}\label{29}
 (\varpi^{\ast}\rho)_{\varepsilon}(v):=(T_{\varepsilon}\varpi)^{-1}[\rho_{\varpi(\varepsilon)}(T_{\varepsilon}\varpi\cdot v)],\;
\forall v \in T_{\varepsilon}^{0}\tilde{G}.
\end{equation}
Finally, for $\sigma$ with the properties stated in the beginning (see (\ref{21})) we consider
\begin{equation}\label{30}
 D(\sigma)=\{(g,\omega)\in\tilde{G}\arrowvert g<\sigma(\pi(g))\;with\; respect\; to\; \omega\} 
\end{equation}
and define for $\psi\in\textsl{C}_{0}^{\infty}\Gamma(\Omega(TG))$ (see (\ref{16}))
\begin{equation}\label{31}
< I(\rho,\sigma,G),\psi>= \int_{ D(\sigma)}\mathrm{d}(\varpi^{\ast}\rho\rightthreetimes\int\psi).
\end{equation}
The density $\mathrm{d}(\varpi^{\ast}\rho\rightthreetimes\int\psi)$ is continuous on $\tilde{G}$, yet it is not plain the convergence 
of the integral on $D(\sigma)$. Remark that the equality $\int(\varphi\circ\pi)\cdot\psi=(\varphi\circ\tilde{\pi})\cdot\int\psi$,  for 
$\varphi\in\textsl{C}_{0}^{\infty}(U)$ and $\psi\in\textsl{C}_{0}^{\infty}\Gamma(\Omega(TG))$, allows to use a partition of unity on $U$
and the Proposition 1 to reduce the analysis to the standard case when $G\subseteq\textbf{R}^{m}\times\textbf{R}$ so that 
$G_{z}=G\cap \pi^{-1}(\{z\})$ be connected $\forall z\in \textbf{R}^{m}$ (possibly void). Here  
$\pi :\textbf{R}^{m}\times\textbf{R}\rightarrow \textbf{R}^{m}$ is the canonical projection.\\
In that case we have for $\rho$ the representation 
\begin{equation}\label{32}
\rho=\lambda\otimes X 
\end{equation}
with $\lambda$ the Lebesgue measure on $\textbf{R}$ and $X$ 
vector field on $G$:
\begin{equation}\label{33}
 X_{(z,y)}=\sum_{i=1}^{m} X^{i}(z,y)\dfrac{\partial}{\partial z_{i}} + X^{m+1}(z,y)\dfrac{\partial}{\partial y},
\end{equation}
$z\in U=\pi(G),\;y\in \textbf{R}$. The section $\sigma$ is then of the form 
\begin{equation}\label{34}
\sigma_{z}=(z,u(z)),
\end{equation}
$(z,u(z))\in G$, where $u:U\rightarrow\textbf{R}$ is Lebesgue measurable and locally essentially bounded. The test density 
\begin{equation}\label{35}
\psi_{(z,y)}=\varphi(z,y)\cdot|\mathrm{d}z_{1}\wedge\dots\wedge\mathrm{d}z_{m}\wedge\mathrm{d}y| 
\end{equation}
is determined by 
$\varphi\in\textsl{C}_{0}^{\infty}(G)$. If we consider $Z^{i}$ such that
\begin{equation}\label{36}
\dfrac{\partial Z^{i}}{\partial y}=X^{i}(z,y),\;1\leqslant i\leqslant m, 
\end{equation}
we get
\begin{multline}\label{37}
< I(\rho,\sigma,G),\psi>=\iint_{G} \mathrm{sgn}(u(z)-y)
\Bigl\{\sum_{i=1}^{m} (Z^{i}(z,u(z))-Z^{i}(z,y))\dfrac{\partial\varphi}{\partial z_{i}}(z,y)+\\
+\Bigl[X^{m+1}(z,u(z))+\sum_{i=1}^{m}\Bigl(\dfrac{\partial Z^{i}}{\partial z_{i}}(z,u(z))-\dfrac{\partial Z^{i}}{\partial z_{i}}(z,y)
\Bigr)\Bigr]\varphi(z,y)\Bigr\}\mathrm{d}z\;\mathrm{d}y.
\end{multline}
Thus the convergence of the integral (\ref{31}) comes from the fact that $\tilde{\pi}(\mathrm{supp}(\int\psi))$ is compact and the section 
$\sigma$ is locally essentially bounded in the sense of (\ref{21}). \\
In the case that the section $\sigma$ is of class $\textsl{C}^{1}$, that is $u$ from (\ref{34}) is $\textsl{C}^{1}$, the generalized
function $I(\rho,\sigma,G)$ is the measurable function given in that standard case by
\begin{equation}\label{38}
I(\rho,\sigma,G)(z,y)= \mathrm{sgn}(u(z)-y)[X^{m+1}(z,u(z))\;-\;\sum_{i=1}^{m}X^{i}(z,u(z))\dfrac{\partial u}{\partial z_{i}}(z)].
\end{equation}
On the other hand, from (\ref{37}) it results that 
$$< I(\rho,\sigma,G),\psi>=\lim_{n\rightarrow\infty}< I(\rho,\sigma_{n},G),\psi>,$$ 
if $\sigma_{n}\tends\limits_{n\to\infty}\sigma$ pointwise almost everywhere and locally uniformly essentially bounded. The 
usual regularization technique allows this type of approximation of a locally essentially bounded function by a sequence of smooth 
functions.\\
For $\zeta\in\textsl{C}_{0}^{\infty}\Gamma(\Omega(TU))$, where $U=\pi(G)$, we consider 
$\overline{\pi^{\ast}}\zeta\in\textsl{C}^{\infty}\Gamma(\Omega(TG/T^{0}G))$ defined by
\begin{equation}\label{39}
 (\overline{\pi^{\ast}}\zeta)_{g}=\zeta_{\pi(g)}\circ(\overline{T_{g}\pi})^{\wedge m}
\end{equation}
where
\begin{equation}\label{40}
\overline{T_{g}\pi}:T_{g}G/T_{g}^{0}G \widetilde{\longrightarrow}T_{\pi(g)}U
\end{equation}
is canonical and $(\overline{T_{g}\pi})^{\wedge m}:(T_{g}G/T_{g}^{0}G)^{\wedge m}\rightarrow(T_{\pi(g)}U)^{\wedge m}$ has also the natural
meaning.\\
Next, for $V=\mathring{V}\subseteq U$ and $\tau$ local $smooth$ section of $\pi$ defined on $V$, with $\tau(V)\subseteq G$, and for
$\zeta\in\textsl{C}_{0}^{\infty}\Gamma(\Omega(TV))$ we define the generalized function on $V$
\begin{equation}\label{41}
< J(\rho,\sigma\arrowvert_{V},\tau),\zeta> =
\int_{\lvert\sigma,\tau\lvert}\mathrm{d}(\rho\rightthreetimes\overline{\pi^{\ast}}\zeta)
+\int_{(\mathrm{im}\;\tau,\;o\;(\tau\rightarrow\sigma))}\rho\rightthreetimes\overline{\pi^{\ast}}\zeta.
\end{equation}
The domain of the first integral is
\begin{equation}\label{42}
 \lvert\sigma,\tau\lvert=:\{g\in G\cap\pi^{-1}(V)\rvert g\in\lvert\sigma_{\pi(g)},\tau_{\pi(g)})\lvert_{G_{\pi(g)}}\}
\end{equation}
(see (\ref{21})). Here $\mathrm{im}\;\sigma$ may pass or jump from one side to the other side of  $\mathrm{im}\;\tau$ but the intervals
$\lvert\sigma_{z},\tau_{z}\lvert_{G_{z}}$ are not oriented; on the other hand the hypersurface integral of an odd $m$-form 
(see (\ref{16})) needs an orientation of the embedding of $\mathrm{im}\;\tau$ in the $m+1$-dimensional domain $G$ (see [5]). The 
symbol $o\;(\tau\rightarrow\sigma)$ denotes the orientation on $G_{z}$ from $\tau_{z}$ to $\sigma_{z}$, if $\sigma_{z}\neq\tau_{z}$, and 
the hypersurface integral is taken only on the region of $\mathrm{im}\;\tau$ where $\sigma_{z}\neq\tau_{z}$. More precisely, if $\omega$ is 
a continuous local orientation of $T^{0}G$, on $W=\mathring{W}\subseteq G$, we consider 
\begin{center}
$s(\tau,\omega,\sigma)(g)=\Biggl\{
\begin{array}{lll}
1,\; \tau_{\pi(g)}<\sigma_{\pi(g)}\\
0,\; \tau_{\pi(g)}=\sigma_{\pi(g)}\\
-1,\;\tau_{\pi(g)}>\sigma_{\pi(g)}
\end{array}$
with \;respect \;to $\omega$ on $G_{\pi(g)}$,
\end{center}
and for $\alpha$ odd $m$-form with support in $W$, we put 
$$\int_{(\mathrm{im}\;\tau,\;o\;(\tau\rightarrow\sigma))}\alpha=\int_{(\mathrm{im}\;\tau,\;\omega)}s(\tau,\omega,\sigma)\cdot \alpha. $$
In the standard framework and with the notations introduced in (\ref{32})-(\ref{34}), for
\begin{equation}\label{43}
 \tau_{z}=(z,v(z)),\; \zeta_{z}=\varphi(z)\cdot |\mathrm{d}z_{1}\wedge\dots\wedge\mathrm{d}z_{m}|,\; z\in V,
\end{equation}
we get
\begin{multline}\label{44}
\int_{\lvert\sigma,\tau\lvert}\mathrm{d}(\rho\rightthreetimes\overline{\pi^{\ast}}\zeta)=\int_{V}\mathrm{sgn}(u(z)-v(z))
\Bigl\{\sum_{i=1}^{m} (Z^{i}(z,u(z))-Z^{i}(z,v(z)))\dfrac{\partial\varphi}{\partial z_{i}}(z)+\\
+\Bigl[X^{m+1}(z,u(z))-X^{m+1}(z,v(z))+\\
+\sum_{i=1}^{m}\Bigl(\dfrac{\partial Z^{i}}{\partial z_{i}}(z,u(z))-
\dfrac{\partial Z^{i}}{\partial z_{i}}(z,v(z))\Bigr)\Bigr]\varphi(z)\Bigr\}\mathrm{d}z,
\end{multline}
\begin{multline}\label{45}
\int_{(\mathrm{im}\;\tau,\;o\;(\tau\rightarrow\sigma))}\rho\rightthreetimes\overline{\pi^{\ast}}\zeta= \int_{V}\mathrm{sgn}(u(z)-v(z))
\Bigl[X^{m+1}(z,v(z))-\\
-\sum_{i=1}^{m} X^{i}(z,v(z))\dfrac{\partial v}{\partial z_{i}}(x)\Big]\varphi(z)\;\mathrm{d}z,
\end{multline}
so that
\begin{multline}\label{46}
 < J(\rho,\sigma\arrowvert_{V},\tau),\zeta>=\int_{V}\mathrm{sgn}(u(z)-v(z))
\Bigl\{\sum_{i=1}^{m} (Z^{i}(z,u(z))-Z^{i}(z,v(z)))\dfrac{\partial\varphi}{\partial z_{i}}(z)+\\
+\Bigl[X^{m+1}(z,u(z))-\sum_{i=1}^{m} X^{i}(z,v(z))\dfrac{\partial v}{\partial z_{i}}(z)+\\
+\sum_{i=1}^{m}\Bigl(\dfrac{\partial Z^{i}}{\partial z_{i}}(z,u(z))-
\dfrac{\partial Z^{i}}{\partial z_{i}}(z,v(z))\Bigr)\Bigr]\varphi(z)\Bigr\}\mathrm{d}z.
\end{multline}
(Recall that $v$ should be smooth - see (\ref{43})).\\
In the case that both $\sigma$ and $\tau$ are smooth, $J(\rho,\sigma\arrowvert_{V},\tau)$ is the measurable function given by
\begin{equation}\label{47}
J(\rho,\sigma\arrowvert_{V},\tau)(z)=\mathrm{sgn}(u(z)-v(z))[X^{m+1}(z,u(z))-
\sum_{i=1}^{m}X^{i}(z,u(z))\dfrac{\partial u}{\partial z_{i}}(z)],
\end{equation}
so that we have (see (\ref{38}))
\begin{equation}\label{48}
J(\rho,\sigma\arrowvert_{V},\tau)= I(\rho,\sigma,G)\circ\tau.
\end{equation}
In the special case when both $\sigma$ and $\tau$ are smooth and defined on $V$ we note the important equalities
\begin{multline}\label{49}
< J(\rho,\sigma,\tau),\zeta>+ < J(\rho,\tau,\sigma),\zeta>=
\int_{\lvert\sigma,\tau\lvert}\mathrm{d}(\rho\rightthreetimes\overline{\pi^{\ast}}\zeta)=\\
=-\int_{(\mathrm{im}\;\tau,\;o\;(\tau\rightarrow\sigma))}\rho\rightthreetimes\overline{\pi^{\ast}}\zeta
-\int_{(\mathrm{im}\;\sigma,\;o\;(\sigma\rightarrow\tau))}\rho\rightthreetimes\overline{\pi^{\ast}}\zeta,
\end{multline}
that come from (\ref{44})-(\ref{47}). For a local section $\tau$ on $V$, with $\tau(V)\subseteq W=\mathring{W}\subseteq G$ 
and $\eta\in\textsl{C}^{\infty}\Gamma(\Omega(TW/T^{0}W))$ we define $\overline{\tau^{\ast}}\eta\in\textsl{C}^{\infty}\Gamma(\Omega(TV))$ by
\begin{equation}\label{50}
 (\overline{\tau^{\ast}}\eta)_{z}=\eta_{\tau_{z}}\circ [(\overline{T_{\tau_{z}}\pi})^{-1}]^{\wedge m}
\end{equation}
(see (\ref{40})).

For the submersion $\pi:F\rightarrow M$, every point from $F$ admits a neighbourhood $W$ and a smooth map $p:W\rightarrow\textbf{R}$ 
such that $(\pi,p):W\rightarrow M\times\textbf{R}$ be a diffeomorphism on an open subset.
Then for $l\in\textsl{C}^{\infty}\Gamma(\Omega(T\textbf{R}))$
it is well defined $\overline{p^{\ast}}l\in\textsl{C}^{\infty}\Gamma(\Omega(T^{0}W))$ by
\begin{equation}\label{51}
 (\overline{p^{\ast}}l)_{g}=l_{p(g)}\circ T_{g}(p\arrowvert_{W_{\pi(g)}}),\; g\in W,
\end{equation}
and if $l_{y}\neq 0,\;\forall y\in\textbf{R}$, then $(\overline{p^{\ast}}l)_{g}\neq 0$, for 
$p\arrowvert_{W_{z}}:W_{z}\rightarrow\textbf{R}$ 
is a diffeomorphism on its open image. Also, for $\psi\in\textsl{C}^{\infty}\Gamma(\Omega(TW)),
\;\; \psi/\overline{p^{\ast}}l\in\textsl{C}^{\infty}\Gamma(\Omega(TW/T^{0}W)$ is defined 
as the factor that satisfies pointwise
\begin{equation}\label{52}
 \psi_{g}=(\overline{p^{\ast}}l)_{g}\rtimes (\psi/\overline{p^{\ast}}l)_{g},\; g\in W.
\end{equation}
The result of disintegration is contained in the following
\begin{theorem} For $W\subseteq G$, $p$ and arbitrary $l$ as before
\begin{equation}\label{53}
< I(\rho,\sigma,G),\psi>=\int_{\textbf{R}}< J(\rho,\sigma,\kappa (y)),
\overline{\kappa (y)^{\ast}}(\psi/\overline{p^{\ast}}l)> \;l(\mathrm{d}y),
\;\;\;\forall\psi\in\textsl{C}_{0}^{\infty}\Gamma(\Omega(TW)),
\end{equation}
where,
\begin{equation}\label{54}
\kappa (y)=:(\pi\arrowvert_{p^{-1}(\{y\})})^{-1},\;\;\forall y\in\textbf{R}.
\end{equation}
Otherwise stated, if $V\times I\subseteq (\pi,p)(W)$ and $\zeta\in\textsl{C}_{0}^{\infty}\Gamma(\Omega(TV))$, then
\begin{equation}\label{55}
< J(\rho,\sigma\arrowvert_{V},\kappa (\cdot)),\zeta>\in \textbf{L}_{loc}^{\infty}(I,\;\textbf{R})
\end{equation}
and
\begin{equation}\label{56}
 \int_{I}< J(\rho,\sigma\arrowvert_{V},\kappa (y)),\zeta>\;l(\mathrm{d}y)=
< I(\rho,\sigma,G),\overline{p^{\ast}}l\rtimes\overline{\pi^{\ast}}\zeta>,
\end{equation}
$\forall\;\zeta\in\textsl{C}_{0}^{\infty}\Gamma(\Omega(TV))$ and $\forall\;l\in \textsl{C}_{0}^{\infty}\Gamma(\Omega(TI))$. Moreover,
if $\{U_{n}\}_{n\in\textbf{N}}$ is a fundamental system of neighbourhoods of $0\in\textbf{R}$ and 
$l_{n}\in\textbf{C}_{0}^{\infty}\Gamma(\Omega(T\textbf{R}))$ with the properties
\begin{center}
$\begin{array}{llllll}
(i) \;\;\mathrm{supp}(l_{n})\in U_{n},\\
\\
(ii)\;\;(l_{n})_{y} \geqslant 0,\;\forall y\in\textbf{R}, \\ 
\\
(iii)\;\;\int_{\textbf{R}} l_{n}(\mathrm{d}y)=1,\\
\\
(iv)\;\;\int_{y<0}l_{n}(\mathrm{d}y)=\int_{y>0}l_{n}(\mathrm{d}y),
\end{array}
\forall n\in\textbf{N},$ \\
\end{center}
then $\forall\zeta\in\textsl{C}_{0}^{\infty}\Gamma(\Omega(TV))$ 
\begin{equation}\label{57}
\lim_{n\rightarrow \infty}\int_{\textbf{R}}< J(\rho,\sigma\arrowvert_{V},\kappa (y+y_{0})),\zeta>\;l_{n}(\mathrm{d}y)= 
< J(\rho,\sigma\arrowvert_{V},\kappa (y_{0})),\zeta>,\;\;\forall y_{0}\in I.
\end{equation}
\end{theorem}

\begin{remark} The property (\ref{57}) does not mean that every point $y_{0}$ is a Lebesgue point for the function (\ref{55}). That 
property is the consequence of the hypotheses on $\{l_{n}\}_{n\in\textbf{N}}$ that ensure also that
\begin{equation}\label{58}
 \lim_{n\rightarrow \infty}\int_{\textbf{R}}\mathrm{sgn}(y)\;l_{n}(\mathrm{d}y)=\mathrm{sgn}(0)=0,
\end{equation}
even if $0$ is not a Lebesgue point for the function $\mathrm{sgn}$.
\end{remark}

In the case of (\ref{32})-(\ref{36}), (\ref{43}), and for $p:W\rightarrow \textbf{R}$ the canonical projection, we have
\begin{multline}\label{59}
< J(\rho,\sigma,\kappa(y)),\zeta>=\int_{V}\mathrm{sgn}(u(z)-y)
\Bigl\{\sum_{i=1}^{m} (Z^{i}(z,u(z))-Z^{i}(z,y))\dfrac{\partial\varphi}{\partial z_{i}}(z)+\\
+\Bigl[X^{m+1}(z,u(z))+\sum_{i=1}^{m}\Bigl(\dfrac{\partial Z^{i}}{\partial z_{i}}(z,u(z))-
\dfrac{\partial Z^{i}}{\partial z_{i}}(z,y)\Bigr)\Bigr]\varphi(z)\Bigr\}\mathrm{d}z 
\end{multline}
(see (\ref{46}) and (\ref{54})). Comparing with (\ref{37}), we see that $I(\rho,\sigma,G)$ is a distribution valued function of $y$, such
that $I(\rho,\sigma,G)(\cdot,y)=J(\rho,\sigma,\kappa(y))$; this is the meaning of the disintegration above in the standard case.\\ 
The composition $I(\rho,\sigma,G)\circ\tau$ from the right hand side of (\ref{48}) has a well defined meaning as soon as $\sigma$ and $\tau$
are locally essentially bounded and \textit{at least one of} $\sigma$ and $\tau$ is \textit{locally Lipschitz continuous}. In that case the 
first identity from (\ref{49})
\begin{equation}\label{60}
< I(\rho,\sigma,G)\circ\tau,\zeta>+ <I(\rho,\tau,G)\circ\sigma,\zeta>=
\int_{\lvert\sigma,\tau\lvert}\mathrm{d}(\rho\rightthreetimes\overline{\pi^{\ast}}\zeta) 
\end{equation}
still holds, as we can see by adding (\ref{46}) with (\ref{47}), the right one of them written for $u$ and $v$ intertwined, and comparing 
the result with (\ref{44}). On the other hand, the generalized function 
\begin{equation}\label{61}
 < C(\rho,\sigma,\tau),\zeta>:=\int_{\lvert\sigma,\tau\lvert}\mathrm{d}(\rho\rightthreetimes\overline{\pi^{\ast}}\zeta)
\end{equation}
is defined for any two essentially locally bounded sections $\sigma$ and $\tau$.
An extension of the above identity to the general case, when $\sigma$ and $\tau$ are only locally essentially bounded, represents the 
key point in the proof of the uniqueness of the entropic solution from Kruzhkov [4]. 
In the standard framework (\ref{32})-(\ref{36}), (\ref{43}), (\ref{59}), it is  considered  
\begin{center}
$\begin{array}{ll}
J(\rho,\sigma,\kappa((p\circ\tau)(w)))(z)+J(\rho,\tau,\kappa((p\circ\sigma) (w)))(z)=\\
=I(\rho,\sigma,G)(z,(p\circ\tau)(w))+I(\rho,\tau,G)(z,(p\circ\sigma)(w)) 
\end{array}$
\end{center}
as a distribution on $V\times V\ni (w,z)$ and the left hand side of (\ref{60}) as the trace on the diagonal of $V\times V$ of this 
distribution; on close inspection the proof of the Theorem 1 from [4] leads to the following
\begin{theorem}
For any locally essentially bounded sections $\sigma$ and $\tau$, in the standard case
\begin{equation}\label{62}
\lim_{h\rightarrow 0}\int< J(\rho,\sigma,\kappa((p\circ\tau)(w)))+J(\rho,\tau,\kappa((p\circ\sigma) (w))),
\zeta(h,w)>\mathrm{d}w=< C(\rho,\sigma,\tau),\zeta>, 
\end{equation}
if
\begin{equation}\label{63}
\zeta(h,w)_{z}=h^{-m}\delta(\dfrac{w-z}{h})\varphi(\dfrac{w+z}{2})\mathrm{d}z
\end{equation}
and $\delta$ is a positive, even and smooth function with compact support and integral equal to 1. Here $\zeta_{z}=\varphi(z)\mathrm{d}z$ as 
in (\ref{43}).
\end{theorem}

\S2.3 \texttt{Entropy density and entropic solutions}

We look at the variational inequality of Kruzhkov [4] as merely a local condition, that a priori contains no distinguished time coordinate 
and no initial datum - as even the author suggests in the last paragraph of the paper [4]. While that inequality reads in [4] as the 
condition $J(\rho,\sigma,\kappa(y))\geqslant 0,\;\forall y$ (according to (\ref{59})), we preferred to consider the equivalent to it 
condition$I(\rho,\sigma,G)\geqslant 0$ concerning the generalized function from (\ref{37}). The Theorem 1 above, apart from an intrinsic 
formulation, allows a more precise statement 
\begin{corollary}
Let $\varrho$, $G$ and $\sigma$ have the properties (\ref{17})-(\ref{21}). If $I(\rho,\sigma,G)\geqslant 0$, on $G$, then 
$J(\rho,\sigma\arrowvert_{V},\tau)\geqslant 0$, on its domain $V$, for an arbitrary smooth local section $\tau$ defined on $V$, with 
$\tau(V)\subseteq G$, and arbitrary $V\subseteq \pi(G)$. Conversely, if $W\subseteq G,\;p:W\rightarrow\textbf{R}$ such that 
$(\pi,p):W\rightarrow M\times\textbf{R}$ be a diffeomorphism on an open subset and 
$J(\rho,\sigma,\kappa (y))\geqslant 0,\;\forall y\in\textbf{R}$ (where $\kappa (y)$ is defined in (\ref{54})), then 
$I(\rho,\sigma,G)\geqslant 0$ on $W$.
\end{corollary}
The following result represents an intrinsic formulation of what, in Theorem 1, \S 3, from [4], appears as key estimate in proving 
the uniqueness of entropic solution submitted to initial conditions. Here it is the consequence of Theorem 2 and Corollary 1:
\begin{corollary}
If $\sigma$ and $\tau$ are locally essentially bounded such that $I(\rho,\sigma,G)\geqslant 0$ and $I(\rho,\tau,G)\geqslant 0$, then 
$C(\rho,\sigma,\tau)\geqslant 0$.
\end{corollary}
While the set of classical solutions of the equation
\begin{equation}\label{64}
\sum_{i=1}^{m}X^{i}(z,u(z))\dfrac{\partial u}{\partial z_{i}}(z)-X^{m+1}(z,u(z))=0 
\end{equation}
(see (\ref{38})), that coincides with the set of $\textsl{C}^{1}$ solutions of the inequation $I(\rho,\sigma,G)\geqslant 0$ , does not 
change when we replace $\rho$ by $f\rho$, for any smooth $f=f(z,y)$ with $f(z,y)\neq 0,\;\forall (z,y)$, the set 
of locally essentially bounded solutions $\sigma$ of the same inequation is invariant only at the multiplication of 
$\rho$ by $f=f(z)$ (i.e. constant on fibers) such that $f(z)>0,\;\forall z$. This fact is better seen when analysing the jump 
admissibility
\begin{proposition}
Let $\varSigma\subset V$ be an embedded surface of class $\textsl{C}^{1}$ that disconnects $V$ in the form 
$V=V_{-}\cup\varSigma\cup V_{+}$ and $\sigma_{-}=:
\sigma\arrowvert_{\overline{V_{-}}}\in\textsl{C}^{1}\Gamma (G\arrowvert_{\overline{V_{-}}}),\;
\sigma_{+}=:\sigma\arrowvert_{\overline{V_{+}}}\in\textsl{C}^{1}\Gamma (G\arrowvert_{\overline{V_{+}}})$. Let 
$\varphi\in\textsl{C}^{1}(V,\textbf{R})$ be such that $V_{-}=\{z\in V\;\arrowvert\;\varphi(z)< 0\},\;V_{+}=
\{z\in V\;\arrowvert\;\varphi(z)>0\},\\
\mathrm{d}_{z}\varphi\neq 0,\;\forall\;z\in \varSigma$. Then $\sigma$ satisfies $I(\rho,\sigma,G)\geqslant 0$ if and only if $\sigma_{-}$
and $\sigma_{+}$ are classical solutions on $V_{-}$ and $V_{+}$ respectively and
\begin{equation}\label{65}
<\mathrm{d}_{z}\varphi,\int_{\lvert k,\sigma_{+}(z)\lvert}T_{g}\pi\cdot \rho^{z}(\mathrm{d}g)>\;\leqslant\;
<\mathrm{d}_{z}\varphi,\int_{\lvert \sigma_{-}(z),k\lvert}T_{g}\pi\cdot \rho^{z}(\mathrm{d}g)>
\end{equation}
$\forall k\in\lvert\sigma_{-}(z),\sigma_{+}(z)\lvert,\;\forall z\in\varSigma;$
here $\rho^{z}=\rho\arrowvert_{G_{z}}$ and the integral is taken on the fiber $G_{z}$ with respect to the  $T_{z}M$-valued vector measure 
$T_{g}\pi\;\rho^{z}(\mathrm{d}g)$.
\end{proposition}
The integral curves of the sub-bundle $D_{g}$ of characteristic directions (see (\ref{18})) are the (reduced) characteristic curves of the 
quasilinear equation (\ref{64}). \textit{The entropy density $\rho$ induces, through its image, which is a half-line in $D_{g}$, an 
orientation on each characteristic curve.} Then also their projections in $M$ are oriented by $\rho$. The result above has then the 
following consequences
\begin{corollary}
In the hypotheses of Proposition 2 we have, $\forall z\in \varSigma$,
the Rankine-Hugoniot condition
\begin{equation}\label{66}
<\mathrm{d}_{z}\varphi, \int_{|\sigma_{+}(z),\sigma_{-}(z)|}T_{g}\pi\;\rho^{z}(\mathrm{d}g)>=0,
\end{equation}
and the relations
\begin{subequations}\label{67}
\begin{eqnarray}
 <\mathrm{d}_{z}\varphi,T_{\sigma_{+}(z)}\pi\;\rho_{\sigma_{+}(z)}(v)>\leqslant 0,\;\forall v\in T_{\sigma_{+}(z)}G_{z},\\
<\mathrm{d}_{z}\varphi,T_{\sigma_{-}(z)}\pi\;\rho_{\sigma_{-}(z)}(v)>\geqslant 0,\;\forall v\in T_{\sigma_{-}(z)}G_{z},
\end{eqnarray}
\end{subequations}
which mean that, moving along the oriented projection of characteristics, one enters, from both sides, in the surface of shock, possibly 
tangentially.
\end{corollary}

\section{General entropy conditions for the equation of 2D flat projective structure}

\S 3.1 \texttt{The extension of the equation to the projective plane}

The equation of 2D flat projective structure is                                                                                                                                                                                                      
\begin{equation}\label{68} 
\dfrac{\partial u}{\partial t}(t,x)+u(t,x)\dfrac{\partial u}{\partial x}(t,x)=0.                                                              
\end{equation}                                                                                                      
The name we use here for this equation is justified first by the following                                                                                                             
\begin{theorem} Let $\Phi$ be a $\textsl{C}^{1}$ diffeomorphism between two open and connected subsets of $\{(t,x,v)\;\arrowvert\;
t,x,v\in\textbf{R}\}$ with the property that maps graphs $\{(t,x,u(t,x))\;\arrowvert\;(t,x)\in D\}$ of local $\textsl{C}^{1}$ 
solutions of (\ref{68}) contained in its domain onto surfaces of the same kind. Then $\exists (s_{j}^{i})_{i,j}$ nonsingular operator in  
$\textbf{R}^{3}$ such that
\begin{equation}\label{69} 
\Phi(t,x,v)=\left(\dfrac{s_{1}^{1}t+s_{2}^{1}x+s_{3}^{1}}{s_{1}^{3}t+s_{2}^{3}x+s_{3}^{3}}\;,\;\dfrac{s_{1}^{2}t+s_{2}^{2}x+s_{3}^{2}}
{s_{1}^{3}t+s_{2}^{3}x+s_{3}^{3}}\;,\;\dfrac{[-\Delta_{3}^{1}t+\Delta_{1}^{1}]v+[\Delta_{3}^{1}x+\Delta_{2}^{1}]}{[-\Delta_{3}^{2}t+
\Delta_{1}^{2}]v+[\Delta_{3}^{2}x+\Delta_{2}^{2}]}\right),                                                                                                      
\end{equation}                                                                                                                                
where $\Delta_{j}^{i}$ is the minor determinant obtained by deleting the $i$-th row and the $j$-th column from the matrix 
$(s_{j}^{i})_{i,j}$. Conversely, every map $\Phi$ of the form  (\ref{69}) sends graphs of local clasical solutions of (\ref{68})
contained in its domain into graphs of local classical solutions.    
\end{theorem}                                                                                                                                 
The point is that if  $(t_{0},x_{0},v_{0})$ lies on the graph of a local classical solution $u$, i.e. $v_{0}=u(t_{0},x_{0})$, then a 
neighbourhood of it, from the (reduced) characteristic  through it, lies on the same graph, i.e. $v_{0}=u(t,x_{0}+v_{0}(t-t_{0})),
\;\forall t$ in a neighbourhood of $t_{0}$. An application $\Phi$ with the property stated in the theorem would map (respective segments
of) characteristics into characteristics and (respective subsets of) fibers of 
$$\textbf{R}^{2}\times\textbf{R}\rightarrow \textbf{R}^{2},\;(t,x,v)\mapsto (t,x),$$ 
into fibers, such that the transformation induced in an open subset of the base  $\textbf{R}^{2}$ should map projections of characteristics 
into projections of characteristics. Our result, surely connected to the M\"{o}bius collineation principle, does not rely on it, however, 
because $\Phi$ is defined only locally with respect to the slope $v$ of the lines from  $\textbf{R}^{2}$. Still, according to the Theorem 3, 
the transformation induced by $\Phi$ in an open subset of the base $\textbf{R}^{2}$ is projective and the action of $\Phi$ on 
$(t_{0},x_{0},v_{0})$ comes from the action of that projective map on the lines $\{(t,x)\;\arrowvert\;x-x_{0}=v_{0}(t-t_{0})\}$ from 
$\textbf{R}^{2}$. Remark that $\Phi$ is \textit{nonlinear in fibers}.\\
This Theorem shows that \textit{the local classical solutions of the quasilinear equation} (\ref{68}) \textit{have to be understood as 
local sections of a nonlinear fiber bundle}. More precisely the right independent of coordinates treatment of the equation (\ref{68})
should proceed as follows: take $P(V)$ the projective plane asociated to a 3-dimensional vector space $V$, i.e.\\ 
$P(V)=\{P(L)\;\arrowvert\;L \;subspace\; of\; V \;of\; 
\dim L=1\}$, and
\begin{equation}\label{70}
G_{2}(V)=:\{P(W)\;\arrowvert\;W subspace\; of\; V\; of \dim W=2\},                                                                 
\end{equation}                                                                                                                             
the Grassmann manifold of the projective lines from $P(V)$. (We may denote also $P(V)=G_{1}(V)$). Next let
\begin{equation}\label{71}
F(V)=:\{(p,d)\;\arrowvert\;p\in P(V),\; d\in G_{2}(V),\;p\in d\} \subset P(V)\times G_{2}(V)),                  
\end{equation}                                                                                                                                
seen as the total space of the fiber bundle
\begin{equation}\label{72} 
\pi:F(V)\longrightarrow P(V),\;\pi:(p,d)\mapsto p,                                                                                     
\end{equation}                                                                                                                                
which plays for the equation of 2D flat projective structure the role of the fiber bundle $\pi:F\rightarrow M$ from the general theory 
presented before. Therefore                                                                                                           
\begin{equation}\label{73}
F(V)_{p_{0}}=\{(p_{0},d)\;\arrowvert\;d\in G_{2}(V), d\ni p_{0}\}.
\end{equation}                                                                                                                                
The characteristic through $(p_{0},d_{0})\in F(V)$ of the equation (\ref{68}) in this setting will be
\begin{equation}\label{74} 
 C_{d_{0}}=\{(p,d_{0})\;\arrowvert\; p\in d_{0}\}.
\end{equation}                                                                                                                                
We may consider also the bundle structure
\begin{equation}\label{75}
 \mathring{\pi}:F(V)\longrightarrow G_{2}(V),\;\mathring{\pi}\;:(p,d)\mapsto d
\end{equation}                                                                                                                                
and remark that $C_{d_{0}}= \mathring{\pi}^{-1}(\{d_{0}\})$.\\
Thus \textit{a local $\textsl{C}^{1}$ section $\sigma_{p}=(p,\delta (p))$ of (\ref{72}) is a solution of our equation if and only if as soon as 
$\delta (p_{0})=d_{0}$ we have also $\delta(p)=d_{0},\;\forall p\in d_{0}$ in a neighbourhood of $p_{0}$}.\\
Coming back to the definition (\ref{74}), we see that the characteristic direction in $(p_{0},d_{0})$ is
\begin{equation}\label{76} 
D_{(p_{0},d_{0})}=: T_{(p_{0},d_{0})}C_{d_{0}}=T_{p_{0}}d_{0}\times \{0_{T_{d_{0}} G_{2}(V)}\}.
\end{equation}                                                                                                                                
We remark now the important one-to-one correspondence
\begin{equation}\label{77} 
\bar{\kappa}_{p_{0}}\;:\{d\in G_{2}(V)\;\arrowvert\; d\ni p_{0}\} \widetilde{\longrightarrow} P(T_{p_{0}}P(V)),\;\bar{\kappa} _{p_{0}}(d)=
P(T_{p_{0}}d),
\end{equation}                                                                                                                                
which establishes also the diffeomorphism
\begin{equation}\label{78}
\kappa_{p_{0}}=:\bar{\kappa}_{p_{0}}\circ\mathring{\pi}\;:F(V)_{p_{0}}\widetilde{\longrightarrow} P(T_{p_{0}}P(V)),
\end{equation}                                                                                                                                
that may be read as a nonlinear fiber bundles (over $P(V)$) isomorphism
\begin{equation}\label{79} 
\kappa\;:F(V)\widetilde{\longrightarrow} P(TP(V)).
\end{equation}                                            
In the general case of an immersion $\iota:C\rightarrow M$ of a curve $C$ in the manifold $M$, its \textit{canonical lift in $P(TM)$} is 
\begin{equation}\label{80}
\Lambda_{\iota}:C\rightarrow P(TM),\;\Lambda_{\iota}(p)=P(T_{p}\iota\;T_{p}C)\in P(T_{\iota(p)}M). 
\end{equation}
Remark that if we identify $F(V)$ with $P(TP(V))$ through $\kappa$, the characteristics of the equation (\ref{68}) are mapped into 
curves like $\kappa(C_{d_{0}})$, which is \textit{the canonical lift of the geodesic} $d_{0}$ in $P(TP(V))$:
\begin{equation}\label{81}
\Lambda_{d_{0}}:d_{0}\rightarrow\kappa(C_{d_{0}}),\;\Lambda_{d_{0}}(p)=P(T_{p}d_{0})=\bar{\kappa}_{p}(d_{0})=\kappa(p,d_{0}). 
\end{equation}
On the other hand, let us consider, for $p_{0}\in P(V),\;A_{(p_{0},\cdot)}$ the tautologic vector bundle over $P(T_{p_{0}}P(V))$: 
for $L\subset T_{p_{0}}P(V)$ subspace of $\dim L=1$
\begin{equation}\label{82}
 A_{(p_{0},P(L))}=:L;
\end{equation}                                                                                                                                
of course $A_{(p_{0},P(L))}$ defines a smooth vector bundle $A$ over $P(TP(V))$. From (\ref{76})-(\ref{79}) we get
\begin{equation}\label{83}
 A_{\kappa (p,d)}=T_{(p,d)}\pi\cdot D_{(p,d)},\;\forall (p,d)\in F(V),
\end{equation}                                                                                                                                
expressing, in particular, the isomorphism between the vector bundle of the characteristic directions in 
the points of a fixed nonlinear fiber $D\arrowvert_{F(V)_{p_{0}}}$ and the tautologic vector bundle $A\arrowvert_{P(T_{p_{0}}P(V))}$.\\

The proof of the Theorem 3 makes appeal to the following \textit{local characterization of a projective transformation through 
its jet of order 2}, that will be useful also later. If $V$ is a real vector space and $S\in GL(V)$ (the group of linear
automorphisms of $V$) we denote $P(S):P(V)\rightarrow P(V)$ the natural map induced in $P(V)$. The set 
\begin{equation}\label{84}
PGL(V)=\{P(S)\arrowvert S\in GL(V)\}
\end{equation}
is the group of projective transformations of $P(V)$. Let $E$ be a real vector space and 
$P(E\times\textbf{R})\supset E$ its canonical projective completion. We say that a function $F:U\rightarrow E$, where $U$ is open in $E$, 
is a \textit{projective map} if $\exists S\in GL(E\times\textbf{R})$ such that $F=P(S)\arrowvert_{U}$. Then the following is true

\textit{$F\in\textsl{C}^{2}(U,E)$ is projective if and only if $\forall z\in U,\;F^{\prime}(z)$ is invertible and 
$\exists \xi(z)\in E^{\ast}$ such that $\forall u,\;v\in E$}
\begin{equation}\label{85}
 F^{\prime}(z)^{-1}\cdot F^{\prime\prime}(z)(u,v)=<\xi(z),u>\cdot v+<\xi(z),v>\cdot u.
\end{equation}
\textit{In that case}
\begin{equation}\label{86}
 F(z+u)=F(z)+\dfrac{1}{1-<\xi(z),u>}\cdot  F^{\prime}(z)u,\;\forall z\in U,\;\forall u\in E \;\arrowvert\; z+u\in U,
\end{equation}
\textit{and $\xi$ verifies}
\begin{equation}\label{87}
\xi(z+u)=\dfrac{1}{1-<\xi(z),u>}\cdot\xi(z),\; \textit{in}\; E^{\ast}, 
\end{equation}
$\forall z,\;z+u\in U$ \textit{with} $1-<\xi(z),u>\neq 0$. \textit{Of course, if (\ref{86}) holds for a fixed $z$ and all 
$u\in E$ with $1-<\xi(z),u>\neq 0$, then $F$ is a projective map}.\\
From (\ref{85}) and (\ref{86}) it results that a projective map is completely determined by its jet of order 2 in a fixed point.

 If we denote as usual $TP(S):TP(V)\rightarrow TP(V)$ the tangent map of $P(S)$ and $P(TP(S)):P(TP(V))\rightarrow P(TP(V))$ its 
projectivized in fibers:
$$P(T_{p}P(S)): P(T_{p}P(V))\rightarrow P(T_{P(S)(p)}P(V)),$$
then the map (see (\ref{77})-(\ref{79}))
$$\kappa^{-1}\circ P(TP(S))\circ \kappa:F(V)\rightarrow F(V)$$
acts naturally on the pairs of a line and a point fixed on it: $(p,d)\mapsto (P(S)(p),P(S)(d))$; here 
$P(S)(d)=\{P(S)(q)\arrowvert q\in d\}$, as usual. 

 Thus the Theorem 3 ensures the existence of $S\in GL(\textbf{R}^{3})$ such that $\Phi=P(TP(S))$, that is, $\Phi$ is the natural lift 
in $P(TP(\textbf{R}^{3}))$ of a projective transformation $P(S)$ of $P(\textbf{R}^{3})$. 

\S 3.2 \texttt{The quasilinear equation of the flat projective structure on a 2D\\
 manifold}

 The question we answer here is the following: how can be characterized the quasilinear equations
\begin{equation}\label{88}
 X^{1}(z,v(z))\dfrac{\partial v}{\partial z_{1}}(z)+
X^{2}(z,v(z))\dfrac{\partial v}{\partial z_{2}}(z)=X^{3}(z,v(z)),\;\;z=(z_{1},z_{2}),
\end{equation}
that are the result of the application of a local bundle diffeomorphism 
\begin{equation}\label{89}
 \Phi(z,y)=(\phi(z),\theta(z,y)),\;\;z=(z_{1},z_{2}),\;\;\phi(z)=(\phi_{1}(z),\phi_{2}(z)), 
\end{equation}
to the equation (see (\ref{68}))
\begin{equation}\label{90}
 \dfrac{\partial u}{\partial z_{1}}(z)+u(z)\dfrac{\partial u}{\partial z_{2}}(z)=0,
\end{equation}
in the sense that 
\begin{equation}\label{91}
v(z)=\theta(\phi^{-1}(z),u(\phi^{-1}(z)))
\end{equation}
satisfies the equation (\ref{88}) as soon as $u$ satisfies (\ref{90})? This question is linked to the fact that the standard Kruzhkov 
entropy condition on the equation (\ref{88}) gives, by transport through $\Phi$ from (\ref{89}), the general entropy condition, in the 
local form, for the equation (\ref{90}).

 We will use the fact that $\Phi$ maps characteristic curves of (\ref{90}) onto characteristic curves of the equation under study 
(\ref{88}). As these characteristic curves are canonical lifts in $P(TM)$ of their projections on the base space $M$, that serve as 
geodesics, for an affine connection without torsion on $M$, say, we have: 1) first to identify and to construct the projective structure 
on a manifold $M$ from these curves in it we call geodesics; 2) next to look for a condition of null curvature for this projective 
structure that would make it locally identical to that of the projective plane; and 3) ensure that the equation (\ref{88}) concerns 
sections of the fiber bundle $P(TU)\rightarrow U$, if $U$ is the open subset of $\textbf{R}^{2}$ where the solutions would be defined, 
that its characteristics are the canonical lifts of their projections in $U$, and that these describe the geodesics of a flat projective 
structure on $U$ according to the results of 1) and 2).

 E. Cartan succeeded to derive the projective structure from its geodesics and to define the suitable projective curvature along with a 
projective connection (see [3]). In fact, the projective structure may be thought as subjacent locally to an affine connection 
without torsion, keeping only the information about the geodesics as immersed curves (and leaving aside the specific affine 
parametrization). Such an affine connection would exist only locally without being unique: two affine connections would be projectively 
equivalent if they define identical geodesics as 1D immersed submanifolds.

 Let $E$ be the vector space that is a local model for the manifold $M$, and $U\subseteq E$ the open image of a local chart on $M$. We 
will denote $\odot$ the symmetric tensor product, such that $(E\odot E)^{\ast}\otimes E$ will stand for the space of symmetric bilinear 
mappings $:E\times E\rightarrow E$. Then an affine connection without torsion on $M$ is described in the domain of the chart by a 
smooth mapping $\Gamma : U\rightarrow (E\odot E)^{\ast}\otimes E$ and the geodesics determined by it are transported by the chart into 
the parametrized curves $\gamma(s)$ in $U$ that satisfy (see [2]):
$$\ddot{\gamma}(s)+\Gamma(\gamma(s))(\dot{\gamma}(s),\dot{\gamma}(s))=0.$$
For a given geodesic with $\dot{\gamma}(s_{0})\neq 0$ we choose $H$ hyperplane in $E$ such that  $\dot{\gamma}(s_{0})\notin H$ and then
fix $e\in E,\;e\notin H$, so that $E=\textbf{R}e\dotplus H$. Let us denote, for $v\in E,\; P^{e}_{H}v\in H$ and
$<P^{H}_{e},v>\in \textbf{R}$ the operators that define the components of $v$ in this decomposition: $v= P^{e}_{H}v+<P^{H}_{e},v>e$. As 
$<P^{H}_{e},\dot{\gamma}(s)>\neq 0,\;\forall s$ in a neighbourhood of $s_{0}$, we may consider the inverse function $s(t)$ such that 
$<P^{H}_{e},\gamma(s(t))>=t,\;\forall t$ in a neighbourhood of $t_{0}=:<P^{H}_{e},\gamma(s_{0})>$. If $g(t)=:P_{H}^{e}\gamma(s(t))\in H$, 
then we have
\begin{equation}\label{92}
 \gamma(s(t))=te+g(t),\;g(t)\in H,
\end{equation}
relation that may serve as definition for both $s(t)$ and $g(t)$. As
\begin{equation}\label{93}
 \dot{s}(t)=\dfrac{1}{<P^{H}_{e},\dot{\gamma}(s(t))>},
\end{equation}
after two derivations of (\ref{92}) we get
\begin{equation}\label{94}
 \ddot{g}(t)=<P^{H}_{e},\Gamma(te+g(t))(e+\dot{g}(t),e+\dot{g}(t))>\dot{g}(t)-P^{e}_{H}\Gamma(te+g(t))(e+\dot{g}(t),e+\dot{g}(t)).
\end{equation}
The equation just found \textit{describes the geodesics as submanifolds}, in fact as graphs of functions (see (\ref{92})). Let us remark 
now that the right hand side of this equality is null if and only if 
$$\Gamma(te+g(t))(e+\dot{g}(t),e+\dot{g}(t))=<P^{H}_{e},\Gamma(te+g(t))(e+\dot{g}(t),e+\dot{g}(t))>(e+\dot{g}(t))$$
which in turn is equivalent to
$$\Gamma(te+g(t))(e+\dot{g}(t),e+\dot{g}(t))\in\textbf{R}(e+\dot{g}(t)).$$
The following statement is true

 \textit{$B\in (E\odot E)^{\ast}\otimes E$ satisfies $B(L\times L)\subseteq L,\;\forall L\in G_{1}(E),$ if and only if 
$\exists \xi\in E^{\ast}$ such that} 
\begin{equation}\label{95}
 B(u,v)=<\xi,u>v+<\xi,v>u,\;\forall u,\;v \in E.
\end{equation}
For a bilinear $B$ defined by $\xi$ as in (\ref{95}) we have $\mathrm{tr}[B(u,\cdot)]=(m+1)<\xi,u>$, if $m=\mathrm{dim}E$; therefore 
the subspace of these symmetric bilinear mappings is the image of the canonical projection in $(E\odot E)^{\ast}\otimes E$
\begin{equation}\label{96}
 [QB](u,v)=(m+1)^{-1}\mathrm{tr}[B(u,\cdot)]\cdot v+(m+1)^{-1}\mathrm{tr}[B(v,\cdot)]\cdot u\;;
\end{equation}
therefore
\begin{equation}\label{97}
 \mathrm{ker} Q=\{B\in(E\odot E)^{\ast}\otimes E\;\arrowvert\;\mathrm{tr}[B(u,\cdot)]=0,\;\forall u\in E\}.
\end{equation}
It results that two local affine connections without torsion $\Gamma_{1}$ and $\Gamma_{2}$ define the same geodesics as immersed 
submanifolds if and only if $(I-Q)(\Gamma_{1}(z)-\Gamma_{2}(z))=0,\;\forall z\in U$; we say in this case that $\Gamma_{1}$ and $\Gamma_{2}$ 
are \textit{projectively equivalent}. And it is natural to define locally, in the respective chart, the \textit{projective connection} by 
$B(z)=:(I-Q)\Gamma(z)$ for any local compatible affine connection without torsion $\Gamma$. If $\chi_{1}$ and $\chi_{2}$ are two charts
defined on the same open subset of $M$ and $B^{\chi_{i}}(\chi_{i}(p))\in\mathrm{ker} Q,\;i=1,\;2,$ are the respectice coefficients defined 
for $p$ in this common domain, then by the simple $(I-Q)$ - projection of the rule for passing from $\Gamma^{\chi_{1}}(\chi_{1}(p))$ to
 $\Gamma^{\chi_{2}}(\chi_{2}(p))$ (see [2]) for an affine connection without torsion $\Gamma$, we get the rule of change for 
the projective connection coefficient
\begin{multline}\label{98}
 B^{\chi_{2}}(\chi_{2}(p))=
(\chi_{2}\circ\chi_{1}^{-1})^{\prime}(\chi_{1}(p))B^{\chi_{1}}(\chi_{1}(p))
[(\chi_{1}\circ\chi_{2}^{-1})^{\prime}(\chi_{2}(p))^{\odot 2}]+\\
+(I-Q)[(\chi_{2}\circ\chi_{1}^{-1})^{\prime}(\chi_{1}(p))(\chi_{1}\circ\chi_{2}^{-1})^{\prime\prime}(\chi_{2}(p))].
\end{multline}
Here we denote, for $B\in(E\odot E)^{\ast}\otimes E$ and $S\in GL(E)$, by $S^{-1} B [S^{\odot 2}]$ 
the pull-back of $B$ through $S$ : $(S^{-1} B [S^{\odot 2}])(u,v)=S^{-1} B(Su,Sv)$ and in short $S^{\ast}B=:S^{-1} B [S^{\odot 2}]$.  
We use the fact that $Q$ \textit{commutes with the pull-back through any} $S\,:\;QS^{\ast}B=S^{\ast}QB$. \\
With this intrinsic projective, but depending 
on the chart, coefficient $B(x)\in \mathrm{ker} Q$, the equation (\ref{94}) of geodesics as immersed curves (\ref{92}) may be written 
\begin{equation}\label{99}
 \ddot{g}(t)=(P^{H}_{e}\otimes\dot{g}(t)-P^{e}_{H})B(te+g(t))(e+\dot{g}(t),e+\dot{g}(t)).
\end{equation}
It is easy to verify next the following statement

 \textit{There is an isomorphism between $\mathrm{ker} Q$ in $E=\textbf{R}e\dotplus H$ and the $H$-valued polynomials of 3-rd degree on 
$H$ of the form}
\begin{equation}\label{100}
 H\ni h\mapsto a(h,h)h+b(h,h)+ch+d\in H,
\end{equation}
\textit{where $a\in (H\odot H)^{\ast},\; b\in(H\odot H)^{\ast}\otimes H,\;c\in H^{\ast}\otimes H,\; d\in H$ are arbitrary coefficients, 
given by}
\begin{equation}\label{101}
 (P^{H}_{e}\otimes h-P^{e}_{H})B(e+h,e+h)=a(h,h)h+b(h,h)+ch+d,\;\;\;\forall h\in H.
\end{equation}
We conclude with the well known \textit{Theorem of E. Cartan: an equation of the form
\begin{equation}\label{102}
 \ddot{g}(t)=a(t,g(t))(\dot{g}(t),\dot{g}(t))\cdot\dot{g}(t)+b(t,g(t))(\dot{g}(t),\dot{g}(t))+c(t,g(t))\dot{g}(t)+d(t,g(t)),
\end{equation}
where $g(t)\in H,\;a(s,h)\in (H\odot H)^{\ast},\,b(s,h)\in (H\odot H)^{\ast}\otimes H,\;c(s,h)\in H^{\ast}\otimes H,\;d(s,h)\in H$, for 
$s\in\textbf{R},\;h\in H$, gives by $t\mapsto te+g(t)$ the geodesics transversal to $H$ in $E=\textbf{R}e\dotplus H$ of an unique 
projective structure on the open set in $E$ where there are defined the local coefficients of the projective connection $a,\;b,\;c,\;d$.}

 For $U$ open subset in $E$ we have $TU=U\times E$ and $P(TU)=U\times P(E)$. On $P(E)$ the splitting $E=\textbf{R}e\dotplus H$ defines a 
standard chart $\chi^{e}_{H}$  
\begin{equation}\label{103}
 (\chi^{e}_{H})^{-1}:H\rightarrow P(E),\;(\chi^{e}_{H})^{-1}(h)=P(\textbf{R}(e+h))\in P(E);
\end{equation}
its domain is $\{P(L)\arrowvert L\in G_{1}(E), L+H=E\}$, the set of transversal to $H$ 1D subspaces of $E$. Then
$\chi^{e}_{H}(P(\textbf{R}(e+\dot{g}(t))))=\dot{g}(t)$, which means that $t\mapsto (te+g(t),\dot{g}(t))$ is the image through
$\mathrm{id}_{U}\times \chi^{e}_{H}$ of the canonical lift in $P(TU)$ of the geodesic $t\mapsto te+g(t)$ from $U$. Taking into account 
that $U$ was the image of a chart $\chi$ on $M$, the equation (\ref{102}) is the transpoted through the chart
\begin{equation}\label{104} 
(\mathrm{id}_{U}\times \chi^{e}_{H})\circ P(T\chi)
\end{equation}
on $P(TM)$ of an equation of the first degree on $P(TM)$ satisfied by the canonical lift of the geodesic. On the other hand it is clear 
that, being given a point in $P(TM)$, hence a point $p\in M$ and a direction $d\in P(T_{p}M)$, there exists a unique geodesic, as 
submanifold, passing through $p$ tangentially to $d$, hence a unique lift of geodesic passing through $(p,d)$. This means that \textit{the 
canonical lifts of geodesics determine a foliation of} $P(TM)$.\\
In order to identify locally on $P(TM)$ this structure we start from a simpler situation. Let $\pi:F\rightarrow M$
be a surjective submersion and $D_{f}\in G_{1}(T_{f}F),\;f\in F$, be a smooth vector sub-bundle of $TF$. The 
immersed curves $C$ in $F$ with $T_{c}C=D_{c},\;\forall c\in C$, are called integral curves of $D$ and the maximal integral curves
are the leaves of the foliation determined by $D$.
Suppose first that
\begin{equation}\label{105}
 T_{f}\pi\;D_{f}\neq 0_{T_{\pi(f)}M},\;\forall f\in F.
\end{equation}
Then \textit{it is well defined the bundle map}
\begin{equation}\label{106}
 \kappa : F\rightarrow P(TM),\; \kappa(f)=P(T_{f}\pi\;D_{f})\in P(T_{\pi(f)}M)
\end{equation}
\textit{on which we impose the condition to be a diffeomorphism on an open subset of} $P(TM)$. This implies, in particular, that
$\mathrm{dim}F=2\mathrm{dim}M-1$. If $\iota:C\rightarrow F$ is an immersed integral curve of $D$, then $\pi\circ\iota$ is an immersion 
and $\kappa\circ\iota$ is the canonical lift of $\pi\circ\iota$. As $\pi$ is a submersion, for every $f_{0}\in F$ there exist a chart of 
$F$ in the neighbourhood of $f_{0}$ and a chart of $M$ in the neighbourhood of $\pi(f_{0})$ that turn $\pi$, under composition on both 
sides, into a linear epimorphism; this means that locally we come to the situation when $F$ is open in 
$\textbf{R}^{m}\times\textbf{R}^{m-1}$, $\pi$ is the restriction to it of the canonical projection
$:\textbf{R}^{m}\times\textbf{R}^{m-1}\rightarrow\textbf{R}^{m}$ and $M$ is the open image of $F$ in $\textbf{R}^{m}$. The properties of $D$
and of $\kappa$ are inherited by restriction to open subsets of $F$ and transferred by diffeomorphisms, so that we may consider them 
fulfilled  also in this standard framework. Here $f=(z,y)$ and 
\begin{equation}\label{107}
 D_{(z,y)}=\textbf{R}\left(\sum_{i=1}^{m}X^{i}(z,y)\dfrac{\partial}{\partial z_{i}}+
\sum_{j=1}^{m-1}X^{j+m}(z,y)\dfrac{\partial}{\partial y_{j}}\right).
\end{equation}
The condition (\ref{105}) means that
\begin{equation}\label{108}
 \sum_{i=1}^{m}X^{i}(z,y)\dfrac{\partial}{\partial z_{i}}\neq 0,\;\forall (z,y)\in F,
\end{equation}
so that, restricting even more our neighbourhood around $f_{0}$ and redefining the chart, we may suppose $X^{1}(z,y)\neq 0,\;\forall 
(z,y)\in F$. In our case the model for the base space is $E=\textbf{R}^{m}$ and in the chart $\chi^{e}_{H}$ on $P(E)$, where $e=e_{1},\;
H=\textbf{R}e_{2}+\ldots +\textbf{R}e_{m}$, we have $\chi^{e}_{H}(P(\textbf{R}(\sum_{i=1}^{m}X^{i}(z,y)\dfrac{\partial}{\partial z_{i}})))=
(X^{2}(z,y)/X^{1}(z,y),\ldots,X^{m}(z,y)/X^{1}(z,y))$. The condition on $\kappa$ becomes: the bundle map
\begin{equation}\label{109}
(z,y)\mapsto (z,(X^{2}(z,y)/X^{1}(z,y),\ldots,X^{m}(z,y)/X^{1}(z,y))) 
\end{equation}
should be a diffeomorphism on its open image. Making this new change of coordinates the sub-bundle $D$ is mapped by $T\kappa$ into a 
1D sub-bundle of $TP(TM)$:
\begin{equation}\label{110}
 (T\kappa\; D)_{(t,x,y)}=\textbf{R}\left(\dfrac{\partial}{\partial t}+\sum_{j=1}^{m-1}y_{j}\dfrac{\partial}{\partial x_{j}}+
\sum_{j=1}^{m-1}Y^{j}(t,x,y)\dfrac{\partial}{\partial y_{j}}\right),
\end{equation}
where we denote $t=z_{1}$, $x_{j}=z_{j+1},\;1\leqslant j\leqslant m-1,\;Y^{j}=:X^{j+m},\;1\leqslant j\leqslant m-1$. So it is defined a map
$Y(t,x,y)$ taking values in $\textbf{R}^{m-1}$, such that the integral curves of the sub-bundle $T\kappa\;D$ in these new coordinates are
of the form $t\mapsto (t,x(t),\dot{x}(t))$, with $x(t)\in\textbf{R}^{m-1}$ and
\begin{equation}\label{111}
 \ddot{x}(t)=Y(t,x(t),\dot{x}(t)).
\end{equation}
Remark that $T\kappa\; D$ is written in (\ref{110}) in the natural coordinates on $P(TM)$ from (\ref{104}), that appear also in (\ref{102}).
Therefore the system defined by $\pi$ and $D$, with the property (\ref{105}) and the condition on $\kappa$ following (\ref{106}), 
corresponds to a projective structure on $M$ if and only if the term $Y$ obtained in (\ref{110}) through the indicated procedure is of the 
form
\begin{equation}\label{112}
 Y(t,x,y)=a(t,x)(y,y)\cdot y+b(t,x)(y,y)+c(t,x)y+d(t,x),
\end{equation}
hence a $\textbf{R}^{m-1}$-valued polynomial of the 3-rd degree in $y\in\textbf{R}^{m-1}$, that is special only in the form of the leading
term, where $a(t,x)(y,y)$ is an arbitrary scalar quadratic form in $y$ (see (\ref{102})).\\
We call the leaves of the foliation above \textit{characteristics of the projective structure}.
The equation just obtained (\ref{102}) determines its coefficients $a,\;b,\;c,\;d$;
it means that \textit{the projective structure is uniquely determined by the foliation of $P(TM)$ through its characteristics}.\\
A smooth local section of $P(TM)\rightarrow M$ is a (classical) solution of the
\textit{quasilinear system of the projective structure} if its image in $P(TM)$ contains, along with any of its points, an entire open 
arc from the unique characteristic through it. Let us write down this system in a local chart on $M$; a smooth solution $\sigma$ is 
transported by the chart $\chi$ in a section of $P(TU)\rightarrow U$, for $U$ open subset in $E$. Let us define $\chi_{\ast}\sigma$ the 
push-forward of $\sigma$ through the chart $\chi$: 
\begin{equation}\label{113}
\chi_{\ast}\sigma=P(T\chi)\circ\sigma\circ\chi^{-1}.
\end{equation}
As $TU=U\times E,\;P(TU)=U\times P(E)$, our section is of the form 
\begin{equation}\label{114}
(\chi_{\ast}\sigma)_{z}=(z,\delta(z)),\;z\in U,\;\delta(z)\in P(E). 
\end{equation}
Restricting the domain of $\chi_{\ast}\sigma$ to the set of points $z$ for which $\delta(z)$ comes from a transversal to $H$ 1D subspace, 
our section will be represented by a $H$-valued function: 
\begin{equation}\label{115}
(\mathrm{id}_{U}\times \chi^{e}_{H})((\chi_{\ast}\sigma)_{z})=(z,u(z)),\;u(z)=\chi^{e}_{H}(\delta(z))\in H. 
\end{equation}
Then the lift of the 
arc of geodesic $\{te+g(t)\arrowvert t\in I\}$ lies on the image of $\sigma$ if and only if 
\begin{equation}\label{116}
\dot{g}(t)=u(te+g(t)),\;\forall t\in I. 
\end{equation}
In the case that a point corresponding to $t_{0}$ of this characteristic belongs to this surface, i.e. $\dot{g}(t_{0})=u(t_{0}e+g(t_{0}))$, 
the equality above is equivalent to its derivative: 
$$\ddot{g}(t)=u^{\prime}(te+g(t))(e+\dot{g}(t)),\;\forall t\in I.$$
Let us denote as in (\ref{102}) $(t,x)=:te+x=z,\;t\in\textbf{R},\;x\in H$, and replace in (\ref{102}) $\dot{g}(t)$ and $\ddot{g}(t)$ 
according to (\ref{116}) and the previous relation; next write the obtained equality for $t_{0}$ and take into account the fact that 
$t_{0},\;g(t_{0}),\;\dot{g}(t_{0})$ were arbitrary. We get the equality
\begin{equation}\label{117}
u^{\prime}(z)(e+u(z))=a(z)(u(z),u(z))\cdot u(z)+b(z)(u(z),u(z))+c(z)u(z)+d(z),
\end{equation}
where $u(z)\in H$; but $u^{\prime}(z)e=\dfrac{\partial u}{\partial t}(t,x) \in H$, $u^{\prime}(z)\arrowvert_{H}=
\dfrac{\partial u}{\partial x}(t,x)\in H^{\ast}\otimes H$, so that finally 
\begin{multline}\label{118}
 \dfrac{\partial u}{\partial t}(t,x)+\dfrac{\partial u}{\partial x}(t,x)u(t,x)=a(t,x)(u(t,x),u(t,x))\cdot u(t,x)+\\
+b(t,x)(u(t,x),u(t,x))+c(t,x)u(t,x)+d(t,x).
\end{multline}
This is \textit{the quasilinear system of the projective structure in a chart}; it is clear that it determines completely this structure
by the coefficients $a,\;b,\;c,\;d$ of the projective connection. The dimension of the system is equal to 
$\mathrm{dim} E-1=\mathrm{dim} H$. Let us analyse in more detail the effect of a change of coordinates, on the base manifold M, for the 
equation (\ref{118}) concerning sections of the fiber bundle $P(TM)\rightarrow M$. If we denote $\phi=:\tilde{\chi}\circ\chi^{-1}$ the 
diffeomorphism, between the open subsets $U$ and $\tilde{U}$ of the model vector space $E$ for $M$, that makes the transition from the 
chart $\chi$ to the chart $\tilde{\chi}$, then 
$$ P(T\phi):U\times P(E)\rightarrow \tilde{U}\times P(E),\;P(T\phi)(z,P(L))=(\phi(z),P(\phi^{\prime}(z)L)),\;z\in U,\;L\in G_{1}(E).$$
Let us consider a splitting $E=\textbf{R}e\dotplus H$. Then (see (\ref{103}))
$$[\chi^{e}_{H}\circ P(\phi^{\prime}(z))\circ (\chi^{e}_{H})^{-1}](h)=\dfrac{P^{e}_{H}\phi^{\prime}(z)e+
P^{e}_{H}\phi^{\prime}(z)h}{<P^{H}_{e},\phi^{\prime}(z)e>+<P^{H}_{e},\phi^{\prime}(z)h>}.$$
If $z=te+x,\;x\in H,\; t\in\textbf{R}$, and $P^{e}_{H}\phi=\xi,\;P^{H}_{e}\phi=\tau$, so that 
$\phi(z)=(\tau(t,x),\xi(t,x))$, we have
$$ [(\mathrm{id}_{\tilde{U}}\times\chi^{e}_{H})\circ P(T\phi)\circ (\mathrm{id}_{U}\times\chi^{e}_{H})^{-1}](z,h)=((\tau(t,x),\xi(t,x)),
\dfrac{\dfrac{\partial\xi}{\partial t}(t,x)+\dfrac{\partial\xi}{\partial x}(t,x)h}{\dfrac{\partial\tau}{\partial t}(t,x)+
\dfrac{\partial\tau}{\partial x}(t,x)h}).$$
For a section $\sigma$ whose expression in the chart $(\mathrm{id}_{U}\times\chi^{e}_{H})\circ P(T\chi)$ is (\ref{115}), the 
expression of $\tilde{u}(\tilde{z})$, for $\tilde{z}=\phi(z)$, in 
the chart $(\mathrm{id}_{U}\times\chi^{e}_{H})\circ P(T\tilde{\chi})$, is obtained according to (see (\ref{113}))
$P(T\phi)\circ\chi_{\ast}\sigma=\tilde{\chi}_{\ast}\sigma\circ\phi$, which means that
\begin{equation}\label{119}
\tilde{u}(\tau(t,x),\xi(t,x))=\dfrac{\dfrac{\partial\xi}{\partial t}(t,x)+
\dfrac{\partial\xi}{\partial x}(t,x)u(t,x)}{\dfrac{\partial\tau}{\partial t}(t,x)+\dfrac{\partial\tau}{\partial x}(t,x)u(t,x)}. 
\end{equation}
Of course the left hand side is defined only there where the denominator from the right hand side is 
$\neq 0$. But $z_{0}\in U$ and $\sigma$ being given, we may choose the splitting 
$E=\textbf{R}e\dotplus H$ such that $\dfrac{\partial\tau}{\partial t}(t_{0},x_{0})+
\dfrac{\partial\tau}{\partial x}(t_{0},x_{0})u(t_{0},x_{0})\neq 0$. Indeed, we take first $e\in E$ such that 
$\delta({z_{0}})=P(\textbf{R}e)$ (see the notation (\ref{114})), next $\alpha\in E^{\ast}$ such that 
$<\alpha,\phi^{\prime}(z_{0})e>\neq 0$ and $<\alpha,e>=1$. Then for $H=:\mathrm{ker}\alpha$ we have 
$P^{H}_{e}=\alpha,\;u(t_{0},x_{0})=0,\; \dfrac{\partial\tau}{\partial t}(t_{0},x_{0})\neq0$.\\ 
Finally the coefficient of the projective connection in the chart $\tilde{\chi}$ is computed from that in the chart $\chi$ according to 
the rule (\ref{98}) and the coefficients $\tilde{a}(\tilde{z}),\;\tilde{b}(\tilde{z}),\;\tilde{c}(\tilde{z}),\;\tilde{d}(\tilde{z})$ are 
obtained using (\ref{101}).\\ 
In the case we are interested here $\mathrm{dim} E=2$, so that $\;\mathrm{dim} H=1$, and the system reduces to one equation; choosing a 
vector as basis in $H$, the coefficients $a,\;b,\;c,\;d$ and $u$ become scalar functions and the equation may be written 
\begin{multline}\label{120}
\dfrac{\partial u}{\partial t}(t,x)+u(t,x)\dfrac{\partial u}{\partial x}(t,x)=a(t,x)u(t,x)^{3}+b(t,x)u(t,x)^{2}+\\
+c(t,x)u(t,x)+d(t,x). 
\end{multline}
Cartan derived a projective curvature differential form (see [3]) that represents, when non zero in a point, the intrinsic 
obstruction to the existence of a chart $\chi$, in a neighbourhood of the point, that would make $B^{\chi}(\chi(p))=0$ for $p$ in the 
domain of the chart (see (\ref{98})). Remark that for such a chart the corresponding coefficients $a,\;b,\;c,\;d$ vanish and the equation 
(\ref{120}) reduces to (\ref{68}). As in the case of the usual curvature of an affine connection, that represents the obstruction to the 
existence of a chart that would make $\Gamma^{\chi}(\chi(p))=0$ on its domain, the existence of the special chart, in the projective case 
also, reduces to a Frobenius complete integrability condition. It results that the projective curvature is null if and only if, around each 
point, there exists a flat affine connection compatible with the projective connection; only that, in the flat projective case, \textit{not 
every} compatible affine connection is flat: the sphere with the usual Levi-Civita connection is projectively flat! As in the affine case, 
where the existence of a null curvature affine connection is equivalent to the existence of an atlas for which the transition 
diffeomorphisms $\chi_{2}\circ\chi_{1}^{-1}$ are all affine maps, in the projective case, the existence of a flat projective structure on a 
manifold is equivalent to the existence of an atlas for which these diffeomorphisms of transition are projective maps (compare 
the characteristic property (\ref{85}), the definition of the projection $Q$ (\ref{96}) and (\ref{98})).\\
The 2D case is somehow special for the computation of that projective curvature form and is not explicit in the book [3];
however it can be done following the lines indicated there. We get the following result

\textit{The 2D projective structure is flat on the domain of the chart if and only if the corresponding to it coefficients
$a,\;b,\;c,\;d$ satisfy on the image of the chart}
\begin{subequations}\label{121}
\begin{eqnarray}
 \dfrac{1}{3}\dfrac{\partial^{2}c}{\partial z_{2}^{2}}-\dfrac{2}{3}\dfrac{\partial^{2}b}{\partial z_{2}\partial z_{1}}+
\frac{\partial^{2}a}{\partial z_{1}^{2}}-\dfrac{2}{3}b\dfrac{\partial b}{\partial z_{1}}+a\dfrac{\partial c}{\partial z_{1}}+
c\dfrac{\partial a}{\partial z_{1}}+\dfrac{1}{3}b\dfrac{\partial c}{\partial z_{2}}-2a\dfrac{\partial d}{\partial z_{2}}
-d\dfrac{\partial a}{\partial z_{2}}=0,\\
\dfrac{1}{3}\frac{\partial^{2}b}{\partial z_{1}^{2}}-\dfrac{2}{3}\dfrac{\partial^{2}c}{\partial z_{1}\partial z_{2}}+
\dfrac{\partial^{2}d}{\partial z_{2}^{2}}
+\dfrac{2}{3}c\dfrac{\partial c}{\partial z_{2}}-d\dfrac{\partial b}{\partial z_{2}}-b\dfrac{\partial d}{\partial z_{2}}-
\dfrac{1}{3}c\dfrac{\partial b}{\partial z_{1}}+2d\dfrac{\partial a}{\partial z_{1}}+a\dfrac{\partial d}{\partial z_{1}}=0;
\end{eqnarray}
\end{subequations}
\textit{(recall that $z_{1}=t,\;z_{2}=x$).} 

 The final result is the following
\begin{theorem}
In order that the equation (\ref{88}) can be transformed through a bundle map (\ref{89}) in the neighbourhood of a point $(z_{0},y_{0})$
into the equation (\ref{120}) it is necessary and sufficient to be fulfilled the following conditions: first (\ref{108}), where $m=2$, 
around $(z_{0},y_{0})$; next, after redefining the indices, the condition (\ref{109}), which in this case means that around $(z_{0},y_{0})$
\begin{equation}\label{122}
 \dfrac{\partial}{\partial y}(X^{2}(z,y)/X^{1}(z,y))\neq 0;
\end{equation}
and finally, should exist the smooth functions $a(z),\;b(z),\;c(z),\;d(z)$ around $z_{0}$ satisfying
\begin{multline}\label{123}
\dfrac{\partial}{\partial z_{1}}(X^{2}(z,y)/X^{1}(z,y))+X^{2}(z,y)/X^{1}(z,y)\dfrac{\partial}{\partial z_{2}}(X^{2}(z,y)/X^{1}(z,y))+\\
+X^{3}(z,y)/X^{1}(z,y)\dfrac{\partial}{\partial y}(X^{2}(z,y)/X^{1}(z,y))=a(z)(X^{2}(z,y)/X^{1}(z,y))^{3}+\\
+b(z)(X^{2}(z,y)/X^{1}(z,y))^{2}+c(z)X^{2}(z,y)/X^{1}(z,y)+d(z),
\end{multline}
$\;\forall (z,y)$ around $(z_{0},y_{0})$. In this case the function
\begin{equation}\label{124}
 u(z)=X^{2}(z,v(z)/X^{1}(z,v(z))
\end{equation}
satisfies the equation (\ref{120}) (where $z_{1}=t,\;z_{2}=x$) as soon as $v$ satisfies the equation (\ref{88}). And the equation
(\ref{120}) can be transformed through a bundle map (\ref{89}) into the equation (\ref{90}) around  $(z_{0},y_{0})$ if and only if
$a(z),\;b(z),\;c(z),\;d(z)$ satisfy the system (\ref{121}) around  $z_{0}$; in that case, there exist a diffeomorphism 
$\phi(z)=(\tau(t,x),\xi(t,x))$, around  $z_{0}$, such that $\tilde{u}$ defined in (\ref{119}) satisfies (\ref{68}) as soon as $u$
satisfies (\ref{120}).
\end{theorem}
\begin{remark}
It may happen that the equation (\ref{88}) can be reduced in the neighbourhood of every point, from the common domain of definition 
of its coefficients, to the equation (\ref{90}), such a transformation, however, being not possible globally. As, for instance, for the 
equation
\begin{equation}\label{125}
 \cos(v(z))\dfrac{\partial v}{\partial z_{1}}(z)+\sin(v(z))\dfrac{\partial v}{\partial z_{2}}(z)=0,
\end{equation}
on the domain 
$$F=\{(z,y)\arrowvert z\in\textbf{R}^{2}, z_{1}>0,\;y\in (-\pi/2,\pi/2)\}\cup\{(z,y)\arrowvert z\in\textbf{R}^{2}, z_{2}>0,
\;y\in (0,\pi)\}.$$
\end{remark}

\S 3.3 \texttt{Restriction to non-linear fiber of entropy densities for \\the equation of 2D flat projective structure}

The equality (\ref{83}) and the remark after (\ref{19}) allow to consider a local entropy density for the equation of 2D flat projective 
structure as a local section $T\pi\;\rho$ of the vector bundle
\begin{equation}\label{126}
 \Omega(T^{0}P(TP(V)))\otimes A\longrightarrow P(TP(V)),                                                                            
\end{equation}                                                                                                                                
where
\begin{equation}\label{127}
(\Omega(T^{0}P(TP(V)))\otimes A)_{(p,q)}=\Omega(T_{q}P(T_{p}P(V)))\otimes A_{(p,q)},\; p\in P(V),\;q\in P(T_{p}P(V)),
\end{equation} 
because the space tangent to the fiber of $F(V)$ is replaced according to 
\begin{equation}\label{128}
T_{(p,d)}\kappa:T^{0}_{(p,d)} F(V)\widetilde {\longrightarrow}T_{\kappa(p,d)}P(T_{p}P(V))
\end{equation}
(see (\ref{77})-(\ref{79}); compare also (\ref{18}), (\ref{76}) and (\ref{83})).\\
It is important to remark that the rule of shock admissibility (\ref{65}) in a point $z$ of the base space $M$ 
involves only the restriction $T\pi\;\rho^{z}$ of the entropy density to the nonlinear fiber $F_{z}$, that appears as 
a local section of the vector bundle $\Omega(TF_{z})\otimes T\pi\;D\arrowvert_{F_{z}}\longrightarrow F_{z}$.                                             
In the case of the equation of 2D flat projective structure $T\pi\;\rho^{p}=:T\pi\;\rho\arrowvert_{P(T_{p}P(V))}$ would be a local 
section of 
$$\Omega(TP(T_{p}P(V)))\otimes A_{(p,\cdot )}\longrightarrow P(T_{p}P(V)).$$                                   
So that we fix, for a while, $p\in P(V)$, denote $E=:T_{p}P(V)$ and consider $A(E)\rightarrow P(E)$,
\begin{equation}\label{129}
 A(E)_{P(L)}=L,\; L\subset E,\; \dim L=1,
\end{equation}
the tautologic vector bundle over $P(E)$. Thus we study local sections $\beta=:T\pi\;\rho^{p}$ (just a notation here) of
\begin{equation}\label{130}
\Omega(TP(E))\otimes A(E)\longrightarrow P(E), 
\end{equation}                                                                                                                                
for $\dim E=2$. Each basis $\{e_{1},e_{2}\}$ in $E$ defines a chart $\chi^{e_{1}}_{e_{2}}$ of the projective line $P(E)$ defined on 
$U_{e_{2}}$, where (see also (\ref{103}))
\begin{equation}\label{131}
U_{e_{2}}=P(E)\smallsetminus \{\textbf{R} e_{2}\};\;\chi^{e_{1}}_{e_{2}}(q)=x,\;q=P(\textbf{R} (e_{1}+x e_{2})),\;x\in\textbf{R}.                                                                   
\end{equation}                                                                                                                                
It will be useful also to consider the local section of $A(E)$:
\begin{equation}\label{132} 
\sigma^{e_{1}}_{e_{2}}(q)=e_{1}+x e_{2}\in L,\;x\in\textbf{R},\;q=P(L)\in U_{e_{2}}.                                                    
\end{equation}                                                                                                                                
Then $\beta\in\textsl{C}\Gamma(\Omega(TP(E))\otimes A(E)\arrowvert_{U})$, where $U\subset U_{e_{2}}$, admits the local 
representation by a scalar function $\beta^{e_{1}}_{e_{2}}$:
\begin{equation}\label{133}
 \beta=(\beta^{e_{1}}_{e_{2}}\circ\chi^{e_{1}}_{e_{2}})\cdot(\chi^{e_{1}}_{e_{2}})^{\ast}\;\lambda\otimes\sigma^{e_{1}}_{e_{2}},
\end{equation}                                                                                                                                
where $\lambda$ stands for the Lebesgue measure on $\textbf{R}$. Remark that $\beta^{e_{1}}_{e_{2}}(x)\neq 0,
\;\forall x\in\chi^{e_{1}}_{e_{2}}(U)$, in virtue of the condition (\ref{17}).\\
If 
\begin{equation}\label{134}
e_{i}^{\prime}=\sum_{j=1}^{2} t_{i}^{j} e_{j},\;1\leqslant i\leqslant 2,
\end{equation}
is a new basis in $E$, then
\begin{equation}\label{135}
 [\chi^{e_{1}}_{e_{2}}\circ (\chi^{e_{1}^{\prime}}_{e_{2}^{\prime}})^{-1}](x^{\prime})=\dfrac{t^{2}_{1}+x^{\prime} t^{2}_{2}}
{t^{1}_{1}+x^{\prime} t^{1}_{2}},
\end{equation}
and finally
\begin{equation}\label{136}
 \beta^{e_{1}^{\prime}}_{e_{2}^{\prime}}(x^{\prime})=\dfrac{|\det T|}{(t^{1}_{1}+x^{\prime} t^{1}_{2})^{3}}\; 
\beta^{e_{1}}_{e_{2}}\left(\dfrac{t^{2}_{1}+x^{\prime} t^{2}_{2}}{t^{1}_{1}+x^{\prime} t^{1}_{2}}\right),
\end{equation}
for $T=(t_{i}^{j})_{i,j}$ from (\ref{134}).\\
As $A(E) \hookrightarrow P(E)\times E$ as a vector sub-bundle over $P(E)$, an entropy density $\beta$ on $U\subset P(E)$ 
defines an $E$-valued measure on $U$. Its push-forward through the chart $\chi^{e_{1}}_{e_{2}}$ will be
\begin{equation}\label{137}
 (\chi^{e_{1}}_{e_{2}})_{\ast}(\beta)_{x}=(e_{1}+x\;e_{2})\beta^{e_{1}}_{e_{2}}(x)\mathrm{d} x.
\end{equation}
For a closed arc $\lvert a,b\rvert\subset U$, with $a\neq b$, it is defined its \textit{barycenter}
\begin{equation}\label{138}
 B_{\beta}(a,b)=:P(\textbf{R}\cdot \beta(\lvert a,b\rvert))\in P(E).
\end{equation}
The barycenter is linked to the Rankine-Hugoniot condition in $p\in \varSigma$ (see (\ref{66})) in the sense that 
$B_{\beta}(\sigma_{+}(p),\sigma_{-}(p))$ determines the tangent direction to the shock curve $\varSigma$ in $p$.
If we denote
\begin{equation}\label{139}
 B^{e_{1}}_{e_{2}}=\chi^{e_{1}}_{e_{2}}\circ B_{\beta}\circ [(\chi^{e_{1}}_{e_{2}})^{-1}\times (\chi^{e_{1}}_{e_{2}})^{-1}],
\end{equation}
from (\ref{137}) we get
\begin{equation}\label{140}
  B^{e_{1}}_{e_{2}}(s,t)=\dfrac{\int_{s}^{t} x \beta^{e_{1}}_{e_{2}}(x)\mathrm{d} x}{\int_{s}^{t} \beta^{e_{1}}_{e_{2}}(x)\mathrm{d} x}.
\end{equation}
As $\beta^{e_{1}}_{e_{2}}$ is continuous and $\beta^{e_{1}}_{e_{2}}(x)\neq 0,\;\forall x\in\chi^{e_{1}}_{e_{2}}(U)$, the function 
$B^{e_{1}}_{e_{2}}$ is of class $\textsl{C}^{1}$ in both arguments, when extended by 
$B^{e_{1}}_{e_{2}}(s,s)=s,\;\forall s\in\chi^{e_{1}}_{e_{2}}(U)$. Moreover $B^{e_{1}}_{e_{2}}$ is strictly increasing, in any of the 
arguments, when the other is fixed; this fact entails the following converse to the Corollary 3, \S 2.3
\begin{proposition}
In the conditions of Proposition 2, \S 2.3, in the case of the equation of 2D flat projective structure, $I(\rho,\sigma,G)\geqslant 0$
if and only if $\sigma_{-}$ and $\sigma_{+}$ are classical solutions on $V_{-}$ and $V_{+}$ respectively, the shock curve $\varSigma$
satisfies the Rankine-Hugoniot condition (\ref{66}), and moving along the oriented (by $\rho$) projection of characteristics, one enters, 
from both sides, in the curve of shock, possibly tangentially (in the sense of (\ref{67})); in fact, for this equation, and 
any fixed $z\in\varSigma$, for any continuous entropy density $\beta=T\pi\;\rho^{z}$ on an open interval of the fiber $P(T_{z}P(V))$ 
including the closed interval $\lvert\sigma_{-}(z),\sigma_{+}(z)\lvert$,  with (\ref{17}) and (\ref{19}), for (\ref{76}) 
(for the fixed $z\in\varSigma$), the Rankine-Hugoniot identity (\ref{66}) and the inequalities (\ref{67}) imply the inequality (\ref{65}), 
$\forall k\in\lvert\sigma_{-}(z),\sigma_{+}(z)\lvert.$
\end{proposition}
In the case that in (\ref{134}) $\textbf{R} e_{2}^{\prime}=\textbf{R} e_{2}=L$, we have identical domains for
$\chi^{e_{1}}_{e_{2}}$ and $\chi^{e_{1}^{\prime}}_{e_{2}^{\prime}}$, i.e. $U_{e_{2}}=U_{e_{2}^{\prime}}$, and in (\ref{134}) 
$t_{2}^{1}=0$; from (\ref{136}) it results that, if $\beta^{e_{1}}_{e_{2}}$ is constant (independent of $x\in\textbf{R}$), then 
$\beta^{e_{1}^{\prime}}_{e_{2}^{\prime}}(x^{\prime})$ is constant again. We call \textit{canonical entropy density on the fiber} 
$P(E)\smallsetminus \{P(L)\}$ a section $\beta$ of (\ref{130}) with this property; it is not unique, but being given a non-zero one, 
any other is a scalar multiple of it. In this case, from (\ref{140}), we get
\begin{equation}\label{141}
  B^{e_{1}}_{e_{2}}(s,t)=\dfrac{s+t}{2},
\end{equation}
and the meaning is that, \textit{for the canonical entropy density $\beta_{P(L)}$ on the fiber $P(E)\smallsetminus \{P(L)\}$ the 
barycenter is determined by the condition that  $B_{\beta_{P(L)}}(a,b)$ is the harmonic conjugate of $P(L)$ with respect to the pair 
$\{a,b\}$ on the projective line $P(E)$}.

\S 3.4 \texttt{The canonical entropy density on the M\"{o}bius band}

 For $V$ real vector space of $\dim V=3$ we fix a point $p_{\infty}\in P(V)$ that will serve as \textit{vanishing point for the 
simultaneity levels} in $P(V)$ and consider
\begin{equation}\label{142}
 M_{p_{\infty}}=:P(V)\smallsetminus\{p_{\infty}\},
\end{equation}
the \textit{M\"{o}bius band}, and $G_{p_{\infty}}\rightarrow M_{p_{\infty}}$ the fiber bundle of \textit{non-instantaneous directions}: for 
$q\in M_{p_{\infty}}$, the strait line $d(q,p_{\infty})$ in $P(V)$ through $q$ and $p_{\infty}$ will be the level of simultaneity of $q,\;
T_{q} d(q,p_{\infty})$ the instantaneous direction from $q$, and
\begin{equation}\label{143}
 (G_{p_{\infty}})_{q}=:P(T_{q} P(V))\setminus\{P(T_{q} d(q,p_{\infty}))\},\;q\in M_{p_{\infty}}.
\end{equation}
For $<\cdot,\cdot>$ scalar product on $V$ and $\omega\in O(V/L_{\infty})$, where $P(L_{\infty})=p_{\infty}$ (see (\ref{4})), we will define 
a smooth entropy density $\rho(p_{\infty}, <\cdot,\cdot>,\omega)$ on $G_{p_{\infty}}$, for the characteristics of the projective structure, 
that induces by restriction to every fiber of $G_{p_{\infty}}$ a canonical entropy density on the fiber. We choose first 
$\{\varepsilon_{1},\varepsilon_{2},\varepsilon_{3}\}$ an orthonormal basis in $V$ with $\textbf{R}\varepsilon_{1}=L_{\infty}$ and 
$\varepsilon_{2}\wedge\varepsilon_{3}$ well-oriented. Next consider the chart on $P(V)$
\begin{equation}\label{144}
 (\kappa^{\varepsilon_{3}}_{\varepsilon_{1},\varepsilon_{2}})^{-1}(x,\theta)=P(\textbf{R}(x\;\varepsilon_{1}+\sin\theta\;\varepsilon_{2}+
\cos\theta\;\varepsilon_{3})),
\end{equation}
where $x\in\textbf{R},\;\theta\in (-\pi/2,\pi/2)$. For the new basis $\varepsilon_{1}^{\prime}=\varepsilon_{1},\;\varepsilon_{2}^{\prime}=
\varepsilon_{3},\;\varepsilon_{3}^{\prime}=-\varepsilon_{2}$, again well-oriented, we have
$$(\kappa^{\varepsilon_{3}^{\prime}}_{\varepsilon_{1}^{\prime},\varepsilon_{2}^{\prime}})^{-1}(x^{\prime},\theta^{\prime})=
(\kappa^{\varepsilon_{3}}_{\varepsilon_{1},\varepsilon_{2}})^{-1}(x,\theta),$$
i.e. $P(\textbf{R}(x^{\prime}\;\varepsilon_{1}+\sin\theta^{\prime}\;\varepsilon_{3}-\cos\theta^{\prime}\;\varepsilon_{2}))=
P(\textbf{R}(x\;\varepsilon_{1}+\sin\theta\;\varepsilon_{2}+\cos\theta\;\varepsilon_{3}))$, only in two cases: 
$(\kappa^{\varepsilon_{3}}_{\varepsilon_{1},\varepsilon_{2}})^{-1}(x,\theta)\in U_{+}$, which means that $x^{\prime}=x,\; \theta^{\prime}=
\theta+\pi/2$, or $(\kappa^{\varepsilon_{3}}_{\varepsilon_{1},\varepsilon_{2}})^{-1}(x,\theta)\in U_{-}$, when 
$x^{\prime}=-x,\;\theta^{\prime}=\theta-\pi/2$. So, on $U_{+}$ we have
$$\dfrac{\partial}{\partial x^{\prime}}=\dfrac{\partial}{\partial x},\;\dfrac{\partial}{\partial \theta^{\prime}}=
\dfrac{\partial}{\partial \theta},$$
while on $U_{-}$
$$\dfrac{\partial}{\partial x^{\prime}}=-\dfrac{\partial}{\partial x},\;\dfrac{\partial}{\partial \theta^{\prime}}=
\dfrac{\partial}{\partial \theta}.$$
For a fixed point $p$ in the domain of $\kappa^{\varepsilon_{3}}_{\varepsilon_{1},\varepsilon_{2}}$ we choose in the tangent plane 
$T_{p}P(V)$ the basis 
\begin{equation}\label{145}
e_{1}=\dfrac{\partial}{\partial \theta},\;e_{2}=\dfrac{\partial}{\partial x},
\end{equation}
and for $p$ in the domain of $\kappa^{\varepsilon_{3}^{\prime}}_{\varepsilon_{1}^{\prime},\varepsilon_{2}^{\prime}}$
the basis
$$e_{1}^{\prime}=\dfrac{\partial}{\partial \theta^{\prime}},\;e_{2}^{\prime}=\dfrac{\partial}{\partial x^{\prime}}.$$
As on both domains $U_{+}$ and $U_{-}$ we have $\textbf{R} e_{2}^{\prime}=\textbf{R} e_{2}$ and $e_{1}^{\prime}=e_{1}$, from (\ref{136})
we infer that, putting $\rho^{e_{1}}_{e_{2}}=1$ on the domain of $\kappa^{\varepsilon_{3}}_{\varepsilon_{1},\varepsilon_{2}}$ and 
$\rho^{e_{1}^{\prime}}_{e_{2}^{\prime}}=1$ on the domain of  
$\kappa^{\varepsilon_{3}^{\prime}}_{\varepsilon_{1}^{\prime},\varepsilon_{2}^{\prime}}$, the definitions will agree on $U_{+}$ and $U_{-}$.
As the domains of these two charts cover $M_{p_{\infty}},\;\rho(p_{\infty}, <\cdot,\cdot>,\omega)$ is well defined, smooth and canonical 
on fibers.\\
Remark that
\begin{equation}\label{146}
 \rho(p_{\infty}, <\cdot,\cdot>,-\omega)=-\rho(p_{\infty}, <\cdot,\cdot>,\omega).
\end{equation}
To see this, take the basis $\{\eta_{1},\eta_{2},\eta_{3}\}$, defined by $\eta_{1}=\varepsilon_{1},\;\eta_{2}=-\varepsilon_{2},\;
\eta_{3}=\varepsilon_{3}$, and note that $(\kappa^{\varepsilon_{3}}_{\varepsilon_{1},\varepsilon_{2}})^{-1}(x,\theta)=
(\kappa^{\eta_{3}}_{\eta_{1},\eta_{2}})^{-1}(x,-\theta)$. Then the choice (\ref{145}) of a basis $\{f_{1},\;f_{2}\}$ in the tangent plane, 
starting from the basis $\{\eta_{1},\eta_{2},\eta_{3}\}$ in $V$, leads to $f_{1}=-e_{1},\;f_{2}=e_{2}$. Using again (\ref{136}), we find 
that $\rho(p_{\infty}, <\cdot,\cdot>,\omega)^{f_{1}}_{f_{2}}=-\rho(p_{\infty}, <\cdot,\cdot>,\omega)^{e_{1}}_{e_{2}}=-1$; on the other hand
$\rho(p_{\infty}, <\cdot,\cdot>,-\omega)^{f_{1}}_{f_{2}}=1$.\\
The orientation $\omega$ of $V/L_{\infty}$ induces an orientation on each subspace $H$ of $V$, supplementary to $L_{\infty}$, hence on 
$P(H)$, since $\dim H=2$; so it induces an orientation on each strait line in $P(V)$ transversal to the simultaneity levels, on each
non-instantaneous direction and finally on each characteristic curve that lies in $G_{p_{\infty}}$. The formula (\ref{137}) and the 
definition in the basis (\ref{145}) shows that this coincides with the orientation on the characteristic curves selected by the
image of $\rho(p_{\infty}, <\cdot,\cdot>,\omega)\arrowvert_{T_{g}P(T_{q}P(V))}$, for $q\in M_{p_{\infty}},\;g\in (G_{p_{\infty}})_{q}$
(see(\ref{143})). Therefore for two scalar products $<\cdot,\cdot>_{1},\;<\cdot,\cdot>_{2}$, there exist $\mu(q)> 0$ such that
$$\rho(p_{\infty}, <\cdot,\cdot>_{2},\omega)\arrowvert_{T_{g}P(T_{q}P(V))}=
\mu(q)\;\rho(p_{\infty}, <\cdot,\cdot>_{1},\omega)\arrowvert_{T_{g}P(T_{q}P(V))},$$ 
or
\begin{equation}\label{147}
 \rho(p_{\infty}, <\cdot,\cdot>_{2},\omega)_{g}=\mu(\pi(g))\;\rho(p_{\infty}, <\cdot,\cdot>_{1},\omega)_{g};
\end{equation}
$\mu$ depends only on $q$ in the base $M_{p_{\infty}}$ since both densities are canonical on the fibers. The remark preceding Proposition 2,
\S 2.3, and the equality (\ref{147}) show that the variational inequality defined by $\rho(p_{\infty}, <\cdot,\cdot>,\omega)$
depends only on $\omega$ and not on the scalar product $<\cdot,\cdot>$ considered.

\S 3.5 \texttt{A simple and meaningful example of entropic solution on the M\"{o}bius\\ band}

We recall first an elementary characteristic property of the conics in the projective plane.\\
\textit{Let $V$ be a real vector space of $\dim V=3,\;b\in (V\odot V)^{\ast}$ a non-degenerate quadratic form with two squares of a sign 
opposite to the sign of the third square, and 
$$C=\{P(\textbf{R} v)\arrowvert v\in V\smallsetminus\{0\},\;b(v,v)=0\},$$
the conic determined by it in the projective plane $P(V)$. Let $p_{1},\;p_{2}\in C,\;p_{1}\neq p_{2},$ and $p_{\infty}$ the intersection 
point of the tangents in $p_{1}$ and $p_{2}$ at $C$. Then for any $q\in C\setminus\{p_{1},p_{2}\}$, the tangent at $C$ in $q$ is the 
harmonic conjugate of $d(q,p_{\infty})$ with respect to the pair $\{d(q,p_{1}),\;d(q,p_{2})\}$.}\\
The following converse is also true.\\
\textit{Let as before $\dim V=3$ and $p_{\infty},\;p_{1},\;p_{2}$ three non-collinear points in $P(V)$. If an arc of curve of class
$\textsl{C}^{1}$ is disjoint from the straight lines $d(p_{\infty},p_{1})$ and $d(p_{\infty},p_{2})$, and for each point $q$ of the arc 
the 
tangent in $q$ at the arc is the harmonic conjugate of the straight line $d(q,p_{\infty})$ with respect to the pair 
$\{d(q,p_{1}),\;d(q,p_{2})\}$, then the arc lies on a conic passing through $p_{1}$, tangentially to $d(p_{\infty},p_{1})$, and through 
$p_{2}$, tangentially to $d(p_{\infty},p_{2})$.}\\
Moreover:\\
\textit{Each point $q\in P(V)\smallsetminus [d(p_{\infty},p_{1})\cup d(p_{\infty},p_{2})\cup d(p_{1},p_{2})]$ lies on a unique conic tangent 
in 
$p_{i}$ to $d(p_{\infty},p_{i})$, for both $i=1,\;2.$ The open arc of this conic of ends $p_{1}$ and $p_{2}$ containing 
$q$ lies completely in the triangle of vertices $p_{1},\;p_{2},\;p_{\infty}$ (there are four such triangles in the projective plane) 
that 
contains $q$ in its interior.}\\
Now our example: let, in the  directly oriented euclidean plane, $p_{\infty}$ be placed in the origin of the axes, $p_{1}$ at both ends 
of the horizontal axis, $p_{2}$ at both ends of the vertical axis; let $r$ be a point in the first quadrant, $p_{1}^{\prime}$ at the 
intersection of $d(r,p_{2})$ with $d(p_{1},p_{\infty})$, and $p_{2}^{\prime}$ at the intersection $d(r,p_{1})\cap d(p_{2},p_{\infty})$.
Next, choose $q_{1}$ on the open segment $(p_{1}^{\prime},r)$ (hence in the first quadrant), and $q_{2}$ on $d(r,p_{1})$ in the second 
quadrant (hence at the left of $p_{2}^{\prime}$). Let us consider the arc of hyperbola from the first quadrant, of asymptotes the two 
axes, passing through $q_{1}$, and similarly the arc of hyperbola from the second quadrant, tangent at $\infty$ to the axes, in 
$p_{1}$ and $p_{2}$ respectively, passing through $q_{2}$. The solution will be defined in the exterior of the quadrilateral 
$p_{2}\;q_{1}\;q_{2}\;p_{1}$ (in direct order around $p_{\infty}$). In the subdomain lying at the right of the straight line 
$d(p_{2},r)$ and under the arc of hyperbola from the first quadrant (itself at the right of $d(p_{2},r)$),we consider a classical 
solution of type ``fan'': the pencil of parallel lines of center $p_{2}$, oriented upwards; in the subdomain bounded by the mentioned 
arc of hyperbola (at the right of $d(p_{2},r)$, in the first quadrant), the segment $[q_{1},q_{2}]$, and the arc of hyperbola from the 
second quadrant lying above the straight line $d(r,p_{1})$, we consider another ``fan'': the pencil of parallel lines of center $p_{1}$, 
oriented from the right to the left; and finally in the subdomain at the left of this arc and above $d(r,p_{1})$, again the pencil of 
parallel lines of center $p_{2}$, but oriented downwards.\\
In order to see better the behaviour of the solution on the straight line from $\infty$, we consider also the case when all the important 
points lie at finite distance. Let $p_{1},\;p_{2},\;r$ be the vertices, ordered by $\omega$, of a triangle in the euclidean plane 
containing $p_{\infty}$ in its interior. Let us denote $p_{1}^{\prime}$ the intersection $d(r,p_{2})\cap d(p_{1},p_{\infty})$ and 
$p_{2}^{\prime}$ at the intersection $d(r,p_{1})\cap d(p_{2},p_{\infty})$. Next, let us consider $q_{1}$ on the open segment 
$(p_{1}^{\prime},r)$ and $q_{2}$ on the open segment $(p_{1},p_{2}^{\prime})$. Consider the triangle of vertices 
$p_{1},\;p_{2},\;p_{\infty}$, in the projective plane, containing $q_{1}$ in its interior, and the arc of conic contained in the interior 
of this triangle passing through $q_{1}$, with ends in $p_{1},\;p_{2}$, tangent to $d(p_{1},p_{\infty})$ and $d(p_{2},p_{\infty})$ 
respectively. We retain from this arc only the segment of ends $q_{1}$ and $p_{1}$; of course, it lies outside the triangle
$p_{1}\;p_{2}\;r$ with $p_{\infty}$ in its interior, and crosses the straight line from $\infty$, so that in the euclidean plane there 
appear only an arc of hyperbola going from $q_{1}$ to $\infty$ and, on the other side of that triangle, an arc from the same hyperbola 
coming from $\infty$ to $p_{1}$. Analogously, consider the triangle of vertices $p_{1},\;p_{2},\;p_{\infty}$, in the projective plane, 
containing $q_{2}$ in its interior, the arc of conic lying in it, passing through $q_{2}$, tangent to $d(p_{1},p_{\infty})$ in $p_{1}$, 
and to $d(p_{2},p_{\infty})$ in $p_{2}$ and retain the segment of conic, from it, of ends $q_{2}$ and $p_{2}$. In the euclidean plane this 
arc, that crosses the straight line from $\infty$, appears as an arc of hyperbola going from $q_{2}$ to $\infty$, and coming from $\infty$,
on the other side of the triangle $p_{1}\;p_{2}\;r$ with $p_{\infty}$ in its interior, on an arc of the same hyperbola, up to $p_{2}$. We 
get in this way, for each vertex of the quadrilateral $p_{2}\;q_{1}\;q_{2}\;p_{1}$, an arc of hyperbola going from the respective vertex
to $\infty$, these disjoint arcs splitting the exterior of this quadrilateral into four disjoint regions, each of them adjacent to the 
quadrilateral by one of its sides. In the region adjacent to the side $[p_{1},p_{2}]$ we consider the pencil of lines of center $p_{1}$,
issuing from $p_{1}$ outside; in the region adjacent to the side $[p_{2},q_{1}]$ the pencil of lines of center $p_{2}$, issuing from $p_{2}$
outside; in the region adjacent to the side $[q_{1},q_{2}]$, again the pencil of lines of center $p_{1}$, but flowing towards $p_{1}$; and 
in the region adjacent by the side $[q_{2},p_{1}]$ the pencil of lines of center $p_{2}$, flowing towards $p_{2}$. The solution is so 
defined in the region of the projective plane exterior to the convex quadrilateral  $p_{2}\;q_{1}\;q_{2}\;p_{1}$ from the euclidean plane, 
containing $p_{\infty}$ in its interior.\\
According to Proposition 3, \S 3.3, and the remark after (\ref{141}), the section of (\ref{72}) constructed above is an entropic 
solution of the equation of 2D flat projective structure with respect to the canonical entropy density 
$$\rho(p_{\infty}, <\cdot,\cdot>,\omega),$$
where $\omega$ is the mentioned orientation around $p_{\infty}$. Remark that, in the case that $q_{1}$ lies on the segment
$(p_{2},p_{1}^{\prime})$ and $q_{2}$ on the segment $(p_{2}^{\prime},r)$, a similar solution with respect to the opposed orientation 
$-\omega$ is obtained.
Its domain is diffeomorphic to the M\"{o}bius band, being the exterior,
in the projective plane, of the contractible neighbourhood of $p_{\infty}$ bounded by the convex quadrilateral 
$p_{2}\;q_{1}\;q_{2}\;p_{1}$ in an affine neighbourhood. Because of that, \textit{in this case}, $\tilde{G}$, covering twice $G$ in 
(\ref{22}), is connected: \textit{the fibers of $G$ are not globally continuously orientable.}

\S 3.6 \texttt{Restriction of entropy densities to non-linear fiber defined by\\ barycentric maps}

\begin{theorem} Let $U\subsetneqq P(E)$ be open and connected, such that the open arc $\lvert p,q\rvert$ is well defined 
$\;\forall p,q\in U,\;p\neq q$. The following nine assertions about a restriction of entropy density to non-linear fiber
$\beta\in\textsl{C}\Gamma(\Omega(TP(E))\otimes A(E)\arrowvert_{U})$ (see (\ref{130})) are equivalent:\\       
$(i)$  $\exists \{e_{1},e_{2}\}$ basis in $E$ with $U\subset U_{e_{2}}$ such that                                                     
\begin{equation}\label{148} 
\int_{\chi^{e_{1}}_{e_{2}}(\lvert p,q\rvert)}\;\beta^{e_{1}}_{e_{2}}(t)\;\texttt{d}t\neq 0,\;\forall p,q\in U,\; p\neq q;                      
\end{equation}                                                                                                                          
$(ii)$  $\forall \{e_{1},e_{2}\}$ basis in $E$ with $U\subset U_{e_{2}}$ property (\ref{148}) holds;\\                                
$(iii)$ $\exists \{e_{1},e_{2}\}$ basis in $E$ with $U\subset U_{e_{2}}$ such that                                                    
\begin{equation}\label{149} 
\forall I\subset\chi^{e_{1}}_{e_{2}}(U),\;\emptyset\neq I=\mathring{I},\;\exists x_{I}\in I\;\arrowvert\; 
\beta^{e_{1}}_{e_{2}}(x_{I})\neq 0,  
\end{equation}                                                                                                                                
and                                                                                                                                           
\begin{equation}\label{150} 
\beta^{e_{1}}_{e_{2}}(x)\;\beta^{e_{1}}_{e_{2}}(y)\geqslant 0,\;\forall x,y\in \chi^{e_{1}}_{e_{2}}(U);                                         
\end{equation}                                                                                                                                
$(iv)$  $\forall \{e_{1},e_{2}\}$ basis in $E$, with $U\subset U_{e_{2}},\;\beta^{e_{1}}_{e_{2}}$ verifies (\ref{149}) and 
(\ref{150});\\ 
$(v)$ 
\begin{equation}\label{151} 
\beta(\lvert p,q\rvert)\neq 0_{E},\;\forall p,q\in U,\;p\neq q, 
\end{equation}
and (see (\ref{138}) for the definition)
\begin{equation}\label{152}
 B_{\beta}(p,q)\in \mathring{\widehat{\lvert p,q\rvert}};
\end{equation}                                                                                                                                
$(vi)$ $\beta$ satisfies (\ref{151}) and $B_{\beta}$ is continuous in the points of $\Delta_{U}=:\{(q,q)\arrowvert q\in U\}$ if
\begin{equation}\label{153}
 B_{\beta}(q,q)=:q,\;\forall q\in U,
\end{equation}                                                                                                                           
(being continuous, in virtue of the hypothesis (\ref{151}), in the points of $U\times U\smallsetminus \Delta_{U}$);\\                          
$(vii)$  $\beta$ satisfies (\ref{151}) and $\forall q\in U$ the function
\begin{equation}\label{154}
 B_{\beta}(q,\cdot)\;:\;U\smallsetminus\{q\}\rightarrow P(E)
\end{equation}                                                                                                                                
is injective;\\                                                                                                                               
$(viii)$  $\exists <\cdot,\cdot>$ scalar product on $E$ such that $\lvert\beta\rvert (\texttt{d}q)$, the total variation 
of the vector measure $\beta$ with respect to the norm defined by the scalar product, verifies
\begin{equation}\label{155} 
\lvert\beta\rvert(J)> 0,\;\forall J\subset U,\;\emptyset\neq J=\mathring{J} 
\end{equation} 
and the function $\omega$, defined $\lvert\beta\rvert$- a. e. by 
$\beta(\texttt{d}q)=\omega(q)\;\lvert\beta\rvert (\texttt{d}q)$, the 
usual Radon-Nikodym decomposition, is continuous on $U$ with values in $E$;\\                                                                       
$(ix)$  $\forall <\cdot,\cdot>$ scalar product on $E$ property (\ref{155}) is fulfilled and  $\omega$ is continuous on $U$.                    
\end{theorem}                                                       
\begin{remark} The characterization of any of the points $(viii)$ and $(ix)$ entails the property $U\neq P(E)$. Indeed, 
$\lVert\omega(q)\rVert=1,\;\forall q\in U$, and $\omega(P(L))\in L,\;\forall P(L)\in U$, as a consequence of the fact that the vector 
measure $\beta$ takes values in the tautologic bundle. But $\nexists\;\omega$ continuous section of the canonical projection of 
$\textbf{S}^{1}$ on $\textbf{P}^{1}$. Or precisely the directions from $P(E)\smallsetminus U$  will play the role of simultaneity 
levels for time evolution.
\end{remark}
In $P(E)$ is defined the \textit{cross-ratio} of four points $[p,q,r,s]$ by
\begin{multline}\label{156}
[P(\textbf{R}\,(e_{1}+x\;e_{2})),P(\textbf{R}\;(e_{1}+e_{2})), P(\textbf{R}\;e_{1}),P(\textbf{R}\;e_{2})]= x,\\
\forall x\in \textbf{R},\;\forall e_{1},\;e_{2}\in E,\;e_{1}\wedge e_{2}\neq 0_{E\wedge E}.                                                                             
\end{multline}                                                                                                                                
Let $\underline{D}\subset P(E)^{4}$ be the maximum domain of definition and continuity of the cross-ratio 
$\underline{c}(a,b,c,d)=
[a,b,c,d],\;\underline{c}:\underline{D}\rightarrow \textbf{R}_{\infty}$, where $\textbf{R}_{\infty}$ is the 
compactification with one point of $\textbf{R}$. Then it is easy to see that $(a,b,c,d)\in\underline{D}$ if and only if the level 
sets of the function $1\mapsto a, 2\mapsto b,3\mapsto c,4\mapsto d$ contain no more than two points. If we denote
\begin{equation}\label{157} 
C_{1}=\underline{D}\smallsetminus\underline{c}^{-1}(\{1\}),\;C_{0}=\underline{D}\smallsetminus\underline{c}^{-1}(\{0\}),\;C_{\infty}=
\underline{D}\smallsetminus\underline{c}^{-1}(\{\infty\})
\end{equation}                                                                                                                                
then $\underline{D}=C_{1}\cup C_{0}=C_{0}\cup C_{\infty}=C_{\infty}\cup C_{1}=(C_{1}\cap C_{0})\cup (C_{0}\cap C_{\infty})
\cup (C_{\infty}\cap C_{1})$.\\                                                                                                                      
We have the following characterization of the barycentric map $B$ associated to a restriction of entropy density to non-linear fiber
of the kind described in the previous theorem:                                                                                                                         
\begin{theorem} Let $E$ be a real vector space of $\dim E=2,\;\emptyset\neq D=\mathring{D}\subset\\
\subset P(E)\times P(E),\;B\;:D\rightarrow P(E)$, be submitted to the following seven hypotheses:\\                                          
$(H1)$ $(p,q)\in D\Rightarrow (q,p)\in D$;\\                                                                                                  
$(H2)$ $(p,q),\; (r,s)\in D\Rightarrow (p,r)\in D$;\\                                                                                         
$(H3)$ $B(p,q)=B(q,p),\;\forall (p,q)\in D$;\\                                                                                             
$(H4)$ $p\neq B(p,q)\neq q,\;\forall (p,q)\in D$ with $p\neq q$;\\
$(H5)$ Let us denote, for $(p,q)\in D$ with $p\neq q$, by $\lvert p,q\rvert$ the closed arc of extremities $p$ and $q$ characterized by
\begin{equation}\label{158}
 B(p,q)\in \mathring{\widehat{\lvert p,q\rvert}}.
\end{equation} 
The hypothesis is that $\forall (p,q)\in D$ with $p\neq q$ and $\forall r,s\in\mathring{\widehat{\lvert p,q\rvert}}$ 
with $r\neq s$ we have $(r,s)\in D$ and $\lvert r,s\rvert\subset\lvert p,q\rvert$;\\
$(H6)$(called the characteristic identity)
\begin{equation}\label{159}
 [B(p,q),B(p,r),B(p,s),B(r,s)]=[B(p,q),B(q,r),B(q,s),B(r,s)],
\end{equation} 
where we suppose both sides being in the case of existence $C_{0}\cap C_{\infty}$ according to (\ref{157});\\
$(H7)$ $B$ continuous on $D$.\\
Then the following consequences hold\\
$(C1)$ $(p,q)\in D\Rightarrow (p,p),\; (q,q)\in D$;\\
$(C2)$ $(q,q)\in D\Rightarrow B(q,q)=q$;\\
$(C3)$ $\exists U\neq\emptyset,\;U=\mathring{U}\subset P(E)$ and $U$ connected such that
\begin{equation}\label{160}
 D=U\times U
\end{equation} 
and
\begin{equation}\label{161}
\lvert p,q\rvert\subset U,\;\forall p,q\in U,\;p\neq q. 
\end{equation} \\
$(C4)$ $\forall p\in U,\; B(p,\cdot)\;:U\rightarrow P(E)$ is injective;\\
$(C5)$ $U\neq P(E)$;\\
$(C6)$ $\forall p\in U,\; B(p,U)\subset U$ and $B(p,\cdot)$ is orientation preserving on $U$;\\
$(C7)$ If we denote
\begin{equation}\label{162}
 r_{q}(a,b,c)=:[B(a,c),B(b,c),B(a,b),q]
\end{equation} 
then $r_{q}$ enjoys of the multiplicative property
\begin{equation}\label{163}
 r_{q}(a,b,d)\;r_{q}(b,c,d)=r_{q}(a,c,d)
\end{equation} 
in the sense that the factors from the left being defined in the case $C_{\infty}$ ($C_{0}$ respectively) it results that 
the right term is defined in the same case $C_{\infty}$ (or $C_{0}$ in the second case);\\
$(C8)$ The difference
\begin{equation}\label{164}
d_{q}(x,y;u,v)=:r_{q}(u,x,v) - r_{q}(u,y,v), 
\end{equation} 
verifies
\begin{equation}\label{165}
 d_{q}(x,y;u,v)\;d_{q}(u,v;w,z)=d_{q}(x,y;w,z)
\end{equation} 
for $ u\neq v,\,w\neq z, \;B(u,v)\neq q\neq B(w,z)$;\\
$(C9)$ For $x\neq y,\; u\neq v,\;\{B(x,y),B(u,v)\}\cap \{p,q\}=\emptyset$ we have
\begin{equation}\label{166}
 \dfrac{d_{p}(x,y;u,v)}{d_{q}(x,y;u,v)}=[B(x,y),B(u,v),p,q].
\end{equation} 
\end{theorem}                                 
                                                                                                                                              
Let us remark that the local section $\sigma^{e_{1}}_{e_{2}}\;: U_{e_{2}}\rightarrow E$ from (\ref{132}) depends only of 
$e_{1}\in E$ and of $\textbf{R}e_{2}\in P(E)$, so that, for $P(L)\in P(E),\;e\in E,\;e\notin L$ we may define 
$U_{P(L)}=: P(E)\smallsetminus\{P(L)\}$ and $\sigma^{e}_{P(L)}\;: U_{P(L)}\rightarrow E$ through its properties
\begin{equation}\label{167}
\sigma^{e}_{P(L_{1})}(P(L_{2}))\in L_{2},\;\sigma^{e}_{P(L_{1})}(P(L_{2}))-e\in L_{1}. 
\end{equation}                                                                                                                                
The following notation will be used: we observe that for $U$ open arc, $U\subsetneqq P(E)$, the right hand side of the relation
\begin{equation}\label{168}
 \epsilon_{U}(p,q;r,s)=:\mathrm{sgn}([p,r,s,k]-[q,r,s,k])
\end{equation}                                                                                                                               
does not depend on $k\in P(E)\smallsetminus U$, if $p,q,r,s\in U$.\\                                                                             
Now we state our result on the determination of the measure $\beta$ through its barycentric map $B$, that will finally lead to the 
formulation of the general variational entropy inequality for the equation of 2D flat projective structure only in terms of the given 
shock rule:                
\begin{theorem} Let $B$ be submitted to $(H1)$-$(H7)$ of the previous theorem and $U$ be defined by its $(C3)$. The equality
\begin{equation}\label{169}
 \beta(\lvert x,y\rvert)=\epsilon_{U}(x,y;u,v)\;d_{q}(x,y;u,v)\;\sigma_{q}^{\beta(\lvert u,v\rvert)}(B(x,y))
\end{equation}                                                                                                                                
for $x\neq y,\;u\neq v,\; B(u,v)\neq q\neq B(x,y)$, may be used first to define $\beta(\lvert x,y\rvert)$, when $u=u_{0},\;v=v_{0}$,
and $\beta(\lvert u,v\rvert)$, for $x=u,\;y=v$ and $u=u_{0},\;v=v_{0}$, by prescribing the values  
$u_{0},\,v_{0},\;\beta(\lvert u_{0},v_{0}\rvert)= \beta_{0}$, submitted to $B(u_{0},v_{0})=P(\textbf{R}\;\beta_{0})$. 
Then the values for $\beta(\lvert x,y\rvert)$ and $\beta(\lvert u,v\rvert)$ found in this way verify the above equality (\ref{169}).
And for this measure we have $B_{\beta}=B$.
\end{theorem}                                                                                                                                 
\begin{remark}
According to the definition (\ref{138}), the barycenter is independent of a non-zero scalar multiple of $\beta$: 
$B_{c\;\beta}=B_{\beta},\;\forall c\neq 0$. In return, the formula (\ref{169}) gives $\beta(\lvert x,y\rvert)$ up to a scalar multiple 
chosen by prescribing $\beta(\lvert u_{0},v_{0}\rvert)=\beta_{0}$, with $P(\textbf{R}\;\beta_{0})=B(u_{0},v_{0})$: we have 
$\sigma^{c\;e}_{q}=c\;\sigma^{e}_{q},\;\forall c\neq 0$.
\end{remark}
The main example of such a measure and corresponding to it barycentric map is that of the already considered canonical entropy density
on the fiber (see the definition before (\ref{141})): let $p\in P(E)$ be fixed; we define $B_{p}$ by
the condition that $B_{p}(x,y)$ should be the harmonic conjugate of $p$ with respect to the pair $\{x,\;y\}$:
\begin{equation}\label{170}
 [B_{p}(x,y),p,x,y]=-1, \forall x\neq p\neq y.
\end{equation}                                                                                                                                
In that case we get
\begin{equation}\label{171}
 [B_{p}(x,y),u,v,p]=\dfrac{[x,u,v,p]+[y,u,v,p]}{2},\; \forall u\neq p\neq v,
\end{equation}                                                                                                                             
$U=U_{p}$ and $\beta$ the Lebesgue measure on $U_{p}$ (which is affine equivalent to $\textbf{R}$).

\end{document}